\documentclass[12pt, twoside]{article}
\usepackage{amsmath,amsthm,amssymb}
\usepackage{times}
\usepackage{enumerate}

\def\titlerunning#1{\gdef\titrun{#1}}
\makeatletter
\def\author#1{\gdef\autrun{\def\and{\unskip, }#1}\gdef\@author{#1}}
\def\address#1{{\def\and{\\\hspace*{18pt}}\renewcommand{\thefootnote}{}%
\footnote {#1}}%
\markboth{\autrun}{\titrun}}
\makeatother
\def\email#1{e-mail: #1}
\def\subjclass#1{{\renewcommand{\thefootnote}{}%
\footnote{\emph{Mathematics Subject Classification (2010):} #1}}}
\def\keywords#1{\par\medskip
\noindent\textbf{Keywords.} #1}

\newtheorem{thm}{Theorem}[section]

\newtheorem{lem}[thm]{Lemma}
\newtheorem{prop}[thm]{Proposition}

\theoremstyle{definition}

\newtheorem{rem}[thm]{Remark}

\numberwithin{equation}{section}

\begin{document}


\baselineskip=17pt


\titlerunning{An observability for parabolic equations
from a measurable set in time}

\title{An observability estimate for parabolic equations \\
from a measurable set in time and its applications}

\author{Kim Dang Phung
\and Gengsheng Wang}

\date{}

\maketitle

\address{K.D. Phung: Math\'{e}matiques -- Analyse, Probabilit\'{e}s, Mod\'{e}lisation -- Orl\'{e}%
ans (MAPMO), Universit\'{e} d'Orl\'{e}ans \& CNRS UMR 6628, F\'{e}d\'{e}%
ration Denis Poisson, Universit\'{e} d'Orl\'{e}ans \& CNRS FR 2964, 45067 Orl%
\'{e}ans Cedex 2, France; \email{kim\_dang\_phung@yahoo.fr} \and G.
Wang: School of Mathematics and Statistics, Wuhan University, Wuhan
430072, China; \email{wanggs62@yeah.net}}



\begin{abstract}
This paper presents a new observability estimate for parabolic
equations in $ \Omega\times\left( 0,T\right) $, where $\Omega$ is a
convex domain. The observation region is restricted over a product
set of an open nonempty subset of $\Omega$ and a subset of positive
measure in $\left( 0,T\right) $. This estimate is derived with the
aid of a quantitative unique continuation at one point in time.
Applications to the bang-bang property for norm and time optimal
control problems are provided.

\end{abstract}

\section{Introduction and main result}

Let $\Omega $ be a bounded, convex and open subset of $\mathbb{R}^{n}$, $%
n\geq 1$, with a boundary $\partial \Omega $. Let $T>0$. We consider
the following parabolic equation:
\begin{equation}
\left\{
\begin{array}{ll}
\partial _{t}u-\Delta u+au+b\cdot \nabla u=0 & \quad \text{in~}\Omega \times
\left( 0,T\right) \ \text{,} \\
u=0 & \quad \text{on~}\partial \Omega \times \left( 0,T\right) \text{ ,} \\
u\left( \cdot ,0\right) \in L^{2}\left( \Omega \right) \text{ .} &
\end{array}%
\right.  \tag{1.1}  \label{1.1}
\end{equation}%
Here $b\in L^{\infty }\left( \Omega \times \left( 0,T\right) \right) ^{n}$, $%
a\in L^{\infty }\left( 0,T;L^{q}\left( \Omega \right) \right) $ with
$q\geq 2 $ for $n=1$, and $q>n$ for $n\geq 2$. Clearly, it defines a
well-posed problem in the sense of Hadamard, that is,
\begin{itemize}
\item for any $u_{0}\in L^{2}\left( \Omega \right) $, there is a unique
solution $u\in C\left( \left[ 0,T\right] ;L^{2}\left( \Omega \right)
\right) $ of (\ref{1.1}) with $u\left( \cdot ,0\right) =u_{0}$;

\item the solution $u$ depends continuously on the initial value.
\end{itemize}
\noindent The above continuous dependence and the uniqueness can be
derived
from the following estimate.%
\begin{equation}
\int_{\Omega }\left\vert u\left( x,t\right) \right\vert ^{2}dx\leq
e^{C_{0}t\left( \left\Vert a\right\Vert _{L^{\infty }\left(
0,T;L^{q}\left( \Omega \right) \right) }^{2}+\left\Vert b\right\Vert
_{L^{\infty }\left( \Omega \times \left( 0,T\right) \right)
}^{2}\right) }\int_{\Omega
}\left\vert u\left( x,0\right) \right\vert ^{2}dx\text{ \quad }\forall t\in %
\left[ 0,T\right] \text{ ,}  \tag{1.2}  \label{1.2}
\end{equation}%
where $C_{0}$ is a positive constant depending only on $\Omega $, $n$ and $q$%
.

This is a kind of stability estimate which shows how the left term $%
\left\Vert u\left( \cdot ,t\right) \right\Vert _{L^{2}\left( \Omega
\right) } $ depends on the right term $\left\Vert u\left( \cdot
,0\right) \right\Vert _{L^{2}\left( \Omega \right) }$. From this
point of view, the
estimate%
\begin{equation}
\left\Vert u\left( \cdot ,T\right) \right\Vert _{L^{2}\left( \Omega
\right) }\leq C_{\left( \Omega ,n,q,\omega ,E,T,a,b\right)
}\int_{D}\left\vert u\left( x,t\right) \right\vert dxdt  \tag{1.3}
\label{1.3}
\end{equation}%
where $D=\omega \times E$ with $\omega $ being an open nonempty subset of $%
\Omega $ and $E$ being a subset of $\left( 0,T\right] $, shows how
the left term $\left\Vert u\left( \cdot ,T\right) \right\Vert
_{L^{2}\left( \Omega \right) }$ depends on the right term
$\left\Vert u\right\Vert _{L^{1}\left( D\right) }$. Here and
throughout the paper, $C_{\left( \cdot \cdot \cdot \right) }$
denotes a positive constant that only depends on what are enclosed
in the brackets. \textit{An interesting problem is to ask what kind
of }$E$\textit{\ makes (\ref{1.3}) standing}.

When $E=\left\{ T\right\} $, (or $E=\left\{ t_{0}\right\} $,
$t_{0}\in \left( 0,T\right] $), (\ref{1.3}) does not hold. However,
it has been obtained (for some potentials $a$ and $b$) that
\begin{equation}
\left\Vert u\left( \cdot ,T\right) \right\Vert _{L^{2}\left( \Omega
\right) }\leq C_{\left( \Omega ,n,q,\omega ,T,a,b\right) }\left\Vert
u\left( \cdot ,0\right) \right\Vert _{L^{2}\left( \Omega \right)
}^{\alpha }\left\Vert u\left( \cdot ,T\right) \right\Vert
_{L^{2}\left( \omega \right) }^{1-\alpha }\text{ ,}  \tag{1.4}
\label{1.4}
\end{equation}%
for some $\alpha \in \left( 0,1\right) $. This is a quantitative
unique continuation at one point in time. It is a kind of H\"{o}lder
continuous dependence in the sense of John. We call (\ref{1.4}) as
the H\"{o}lder continuous dependence from one point in time. With
regard to the studies of unique continuation, we refer the readers
to \cite{BT}, \cite{L}, \cite{K}, \cite{KT} and references therein.

When $E=\left( 0,T\right) $ (or $E$ is a subinterval of $\left( 0,T\right) $%
), the estimate (\ref{1.3}), viewed as a refined observability
estimate in control theory of PDE, has been discussed in many
literatures (see for
instance \cite{LR}, \cite{FI}, \cite{DZZ}). It is obtained that the estimate (%
\ref{1.3}) holds for a large class of potentials $a$ and $b$ (see \cite{DFGZ}%
).

The present paper studies the estimate (\ref{1.3}) when $E$ is a
measurable set of $\left( 0,T\right) $ with a positive measure. The
main result is presented as follows.

\begin{thm} \label{theorem1.1} Let $E\subset \left( 0,T\right) $\ be a measurable set
with a positive measure. Let $\omega $\ be a nonempty open subset of
$\Omega $. Then any solution $u$ to (\ref{1.1}) holds the estimate
\begin{equation}
\left\Vert u\left( \cdot ,T\right) \right\Vert _{L^{2}\left( \Omega
\right) }\leq C_{\left( \Omega ,n,q,\omega ,E,T,a,b\right)
}\int_{\omega \times E}\left\vert u\left( x,t\right) \right\vert
dxdt\text{ .}  \tag{1.5} \label{1.5}
\end{equation}

\end{thm}

The key to establish Theorem \ref{theorem1.1} is the following
strategy:
\begin{align*}
& \text{H\"{o}lder continuous dependence from one point in time } \\
& \Longrightarrow \text{ Observability from a measurable set in time
\quad (i.e.,(\ref{1.5})).}
\end{align*}%
This method allows us to build up (\ref{1.5}) for parabolic
equations with space-time dependent potentials $a$ and $b$. It also
provides a different
way from that in \cite{W} to get (\ref{1.5}) for the case where $a=0$, $b=0$%
. The above-mentioned strategy is partially inspired by \cite{M}. In
our paper, the estimate (\ref{1.4}) is built up by the technique
provided in \cite{P}, \cite{EFV} and \cite{PW}.

The rest of the paper is organized as follows. Section 2 first shows the H%
\"{o}lder continuous dependence from one point in time, and then
presents the proof of Theorem \ref{theorem1.1}. Section 3 provides
some applications of Theorem \ref{theorem1.1} to the bang-bang
property for norm and time optimal control problems. In Appendix,
the proof of some results (which are used in the proof of Theorem
\ref{theorem1.1}) is given.

\section{Proof of Theorem \ref{theorem1.1}}

\subsection{Preliminary results}

The proof of Theorem \ref{theorem1.1} is based on the following two
results. We provide the proof of the first one in Appendix and that
of the second one in subsection 2.3.

\begin{prop} \label{proposition2.1}Let $E\subset \left( 0,T\right) $\ be a measurable
set with a positive measure. Let $\ell $\ be a density point for
$E\subset \left( 0,T\right) $. Then for each $z>1$, there exists a
$\ell _{1}\in \left( \ell ,T\right) $\ such that the sequence
$\left\{ \ell _{m}\right\} _{m\geq 1}$, given by
\begin{equation}
\ell _{m+1}=\ell +\frac{1}{z^{m}}\left( \ell _{1}-\ell \right)
\text{ ,} \tag{2.1.1}  \label{2.1.1}
\end{equation}%
satisfies%
\begin{equation}
\ell _{m}-\ell _{m+1}\leq 3\left\vert E\cap \left( \ell _{m+1},\ell
_{m}\right) \right\vert \text{ .}  \tag{2.1.2}  \label{2.1.2}
\end{equation}
\end{prop}

To state the second result, we need the following notation. Let
\begin{equation*}
p=\left\vert
\begin{array}{ll}
\frac{2n}{q} & \text{if }n<q\leq 2n \\
1 & \text{if }2n\leq q\text{ .}%
\end{array}%
\right.
\end{equation*}%
Write
\begin{equation*}
A\left( T,\left\Vert a\right\Vert \right) =\left\Vert a\right\Vert
_{L^{\infty }\left( 0,T;L^{q}\left( \Omega \right) \right) }+\left(
T+T^{2-p}\right) \left\Vert a\right\Vert _{L^{\infty }\left(
0,T;L^{q}\left( \Omega \right) \right) }^{2}+T^{2}\left( \left\Vert
a\right\Vert _{L^{\infty
}\left( 0,T;L^{q}\left( \Omega \right) \right) }\right) ^{\frac{4}{2-p}}%
\text{,}
\end{equation*}%
\begin{equation*}
K\left( T,\left\Vert a\right\Vert ,\left\Vert b\right\Vert \right)
=1+A\left( T,\left\Vert a\right\Vert \right) +T\left\Vert
b\right\Vert _{L^{\infty }\left( \Omega \times \left( 0,T\right)
\right) }^{2}\text{ }
\end{equation*}%
and
\begin{equation*}
\beta \left( r,T,\left\Vert b\right\Vert \right)
=\frac{1}{r^{2}}e^{2T\left( 1+\left\Vert b\right\Vert _{L^{\infty
}\left( \Omega \times \left( 0,T\right) \right) }^{2}\right) }\text{
.}
\end{equation*}

\begin{prop} \label{proposition2.2}Let $B_{r}$\ be an open ball of radius $r>0$\ and
contained in $\Omega $. There is a $C=C_{(\Omega ,n,q)}$\ such that
any solution $u$\
to (\ref{1.1}) satisfies%
\begin{equation}
\begin{array}{ll}
\displaystyle\int_{\Omega }\left\vert u\left( x,L\right) \right\vert
^{2}dx & \leq \left( C\displaystyle\int_{B_{r}}\left\vert u\left(
x,L\right) \right\vert ^{2}dx\right) ^{1-\alpha \left(
r,T,\left\Vert b\right\Vert
\right) } \\
& \quad \times \left( e^{C\left( K\left( T,\left\Vert a\right\Vert
,\left\Vert b\right\Vert \right) +\frac{1}{L}\right) }\displaystyle%
\int_{\Omega }\left\vert u\left( x,0\right) \right\vert
^{2}dx\right)
^{\alpha \left( r,T,\left\Vert b\right\Vert \right) }%
\end{array}
\tag{2.1.3}  \label{2.1.3}
\end{equation}%
where $L$\ is arbitrarily taken from $\left( 0,T\right] $, and where%
\begin{equation*}
\alpha \left( r,T,\left\Vert b\right\Vert \right) =\frac{C\beta
\left( r,T,\left\Vert b\right\Vert \right) }{1+C\beta \left(
r,T,\left\Vert b\right\Vert \right) }\text{ .}
\end{equation*}%
Furthermore, there is a positive constant $c$\ (only depending on
$\Omega $, $n$\ and $q$) such that any solution $u$\ to (\ref{1.1})
satisfies
\begin{equation}
\begin{array}{ll}
\left\Vert u\left( \cdot ,t_{2}\right) \right\Vert _{L^{2}\left(
\Omega \right) } & \leq \frac{1}{\varepsilon ^{\gamma \left(
r,T,\left\Vert b\right\Vert \right) }}e^{c\left( K\left(
T,\left\Vert a\right\Vert ,\left\Vert b\right\Vert \right)
+\frac{1}{t_{2}-t_{1}}\right) \beta \left( r,T,\left\Vert
b\right\Vert \right) }\left\Vert u\left( \cdot ,t_{2}\right)
\right\Vert _{L^{1}\left( B_{r}\right) } \\
& \quad +\varepsilon \left\Vert u\left( \cdot ,t_{1}\right)
\right\Vert
_{L^{2}\left( \Omega \right) }\quad \forall \varepsilon >0%
\end{array}
\tag{2.1.4}  \label{2.1.4}
\end{equation}%
where $t_{1}$\ and $t_{2}$\ are arbitrarily taken such that $0\leq
t_{1}<t_{2}\leq T$, and where%
\begin{equation}
\gamma \left( r,T,\left\Vert b\right\Vert \right) =C\beta \left(
r/2,T,\left\Vert b\right\Vert \right) \left( 1+n/2\right) +n/2\text{
.} \tag{2.1.5}  \label{2.1.5}
\end{equation}
\end{prop}

\subsection{Proof of Theorem \ref{theorem1.1}}

Write $B_{r}$ for an open ball of radius $r>0$ and contained in
$\omega $.
Let $\ell $ be a density point for $E\subset \left( 0,T\right) $. Let $%
\left\{ \ell _{m}\right\} _{m\geq 1}$ be the sequence provided by
Proposition \ref{proposition2.1} with $z=\sqrt{\frac{\gamma
+2}{\gamma +1}}$, where $\gamma $ is given by (\ref{2.1.5}). Let
$t\in \left( \ell _{m+1},\ell _{m}\right] $.
Then we apply (\ref{2.1.4}) in Proposition \ref{proposition2.2}, where $t_{2}=t$ and $%
t_{1}=\ell _{m+2}$, to get that
\begin{equation}
\begin{array}{ll}
\left\Vert u\left( \cdot ,t\right) \right\Vert _{L^{2}\left( \Omega
\right) } & \leq \frac{1}{\varepsilon ^{\gamma }}e^{c\left( K\left(
T,\left\Vert
a\right\Vert ,\left\Vert b\right\Vert \right) +\frac{1}{t-\ell _{m+2}}%
\right) \beta \left( r,T,\left\Vert b\right\Vert \right) }\left\Vert
u\left(
\cdot ,t\right) \right\Vert _{L^{1}\left( B_{r}\right) } \\
& \quad +\varepsilon \left\Vert u\left( \cdot ,\ell _{m+2}\right)
\right\Vert _{L^{2}\left( \Omega \right) }\quad \forall \varepsilon
>0\text{
.}%
\end{array}
\tag{2.2.1}  \label{2.2.1}
\end{equation}%
Since it follows from (\ref{1.2}) that
\begin{equation*}
\left\Vert u\left( \cdot ,\ell _{m}\right) \right\Vert _{L^{2}\left(
\Omega \right) }\leq e^{C_{0}T\left[ \left\Vert a\right\Vert
_{L^{\infty }\left( 0,T,L^{q}\left( \Omega \right) \right)
}^{2}+\left\Vert b\right\Vert _{L^{\infty }\left( \Omega \times
\left( 0,T\right) \right) }^{2}\right] }\left\Vert u\left( \cdot
,t\right) \right\Vert _{L^{2}\left( \Omega \right) }\text{ ,}
\end{equation*}%
we integrate (\ref{2.2.1}) over $E\cap \left( \ell _{m+1},\ell
_{m}\right) $
to get that%
\begin{equation*}
\begin{array}{ll}
& \quad \left\vert E\cap \left( \ell _{m+1},\ell _{m}\right)
\right\vert e^{-C_{0}T\left[ \left\Vert a\right\Vert _{L^{\infty
}\left( 0,T,L^{q}\left( \Omega \right) \right) }^{2}+\left\Vert
b\right\Vert _{L^{\infty }\left( \Omega \times \left( 0,T\right)
\right) }^{2}\right] }\left\Vert u\left(
\cdot ,\ell _{m}\right) \right\Vert _{L^{2}\left( \Omega \right) } \\
& \leq \frac{1}{\varepsilon ^{\gamma }}e^{c\left( K\left(
T,\left\Vert a\right\Vert ,\left\Vert b\right\Vert \right)
+\frac{1}{\ell _{m+1}-\ell
_{m+2}}\right) \beta \left( r,T,\left\Vert b\right\Vert \right) }%
\displaystyle\int_{E\cap \left( \ell _{m+1},\ell _{m}\right)
}\left\Vert
u\left( \cdot ,t\right) \right\Vert _{L^{1}\left( B_{r}\right) }dt \\
& \quad +\left\vert E\cap \left( \ell _{m+1},\ell _{m}\right)
\right\vert \varepsilon \left\Vert u\left( \cdot ,\ell _{m+2}\right)
\right\Vert
_{L^{2}\left( \Omega \right) }\quad \forall \varepsilon >0\text{ .}%
\end{array}%
\end{equation*}%
This, along with (\ref{2.1.1}) and (\ref{2.1.2}), indicates that
there is a
positive constant $d$ (only depending on $\Omega $, $n$ and $q$) such that%
\begin{equation*}
\begin{array}{ll}
& \quad \left\Vert u\left( \cdot ,\ell _{m}\right) \right\Vert
_{L^{2}\left(
\Omega \right) } \\
& \leq \frac{1}{\varepsilon ^{\gamma }}e^{dK\left( T,\left\Vert
a\right\Vert ,\left\Vert b\right\Vert \right) \beta \left(
r,T,\left\Vert b\right\Vert \right) }\left[ \frac{1}{\ell _{1}-\ell
}\frac{z^{m}}{z-1}\right] e^{d\beta
\left( r,T,\left\Vert b\right\Vert \right) \left[ \frac{1}{\ell _{1}-\ell }%
\frac{z^{m+1}}{z-1}\right] } \\
& \quad \times \displaystyle\int_{E\cap \left( \ell _{m+1},\ell
_{m}\right) }\left\Vert u\left( \cdot ,t\right) \right\Vert
_{L^{1}\left( B_{r}\right) }dt+\varepsilon \left\Vert u\left( \cdot
,\ell _{m+2}\right) \right\Vert
_{L^{2}\left( \Omega \right) } \\
& \leq \frac{1}{\varepsilon ^{\gamma }}e^{dK\left( T,\left\Vert
a\right\Vert ,\left\Vert b\right\Vert \right) \beta \left(
r,T,\left\Vert b\right\Vert \right) }e^{\left( 1+d\beta \left(
r,T,\left\Vert b\right\Vert \right)
\right) \left[ \frac{1}{\ell _{1}-\ell }\frac{z^{m+1}}{z-1}\right] }%
\displaystyle\int_{E\cap \left( \ell _{m+1},\ell _{m}\right)
}\left\Vert
u\left( \cdot ,t\right) \right\Vert _{L^{1}\left( B_{r}\right) }dt \\
& \quad +\varepsilon \left\Vert u\left( \cdot ,\ell _{m+2}\right)
\right\Vert _{L^{2}\left( \Omega \right) }\quad \forall \varepsilon
>0\text{
,}%
\end{array}%
\end{equation*}%
that is
\begin{equation}
\begin{array}{ll}
& \quad \varepsilon ^{\gamma }e^{-\eta z^{m+2}}\left\Vert u\left(
\cdot ,\ell _{m}\right) \right\Vert _{L^{2}\left( \Omega \right)
}-\varepsilon ^{\gamma +1}e^{-\eta z^{m+2}}\left\Vert u\left( \cdot
,\ell _{m+2}\right)
\right\Vert _{L^{2}\left( \Omega \right) } \\
& \leq e^{dK\left( T,\left\Vert a\right\Vert ,\left\Vert
b\right\Vert
\right) \beta \left( r,T,\left\Vert b\right\Vert \right) }\displaystyle%
\int_{E\cap \left( \ell _{m+1},\ell _{m}\right) }\left\Vert u\left(
\cdot ,t\right) \right\Vert _{L^{1}\left( B_{r}\right) }dt\quad
\forall
\varepsilon >0\text{ ,}%
\end{array}
\tag{2.2.2}  \label{2.2.2}
\end{equation}%
where $\eta =\left( 1+d\beta \left( r,T,\left\Vert b\right\Vert
\right) \right) \left[ \frac{1}{\ell _{1}-\ell }\frac{1}{z\left(
z-1\right) }\right] $. By taking $\varepsilon =e^{-\eta z^{m+2}}$ in
(\ref{2.2.2}), and by using the fact that $\left( \gamma +1\right)
z^{2}=\gamma +2$, we obtain that
\begin{equation}
\begin{array}{ll}
& \quad e^{-\eta \left( \gamma +2\right) z^{m}}\left\Vert u\left(
\cdot ,\ell _{m}\right) \right\Vert _{L^{2}\left( \Omega \right)
}-e^{-\eta \left( \gamma +2\right) z^{m+2}}\left\Vert u\left( \cdot
,\ell _{m+2}\right)
\right\Vert _{L^{2}\left( \Omega \right) } \\
& \leq e^{dK\left( T,\left\Vert a\right\Vert ,\left\Vert
b\right\Vert
\right) \beta \left( r,T,\left\Vert b\right\Vert \right) }\displaystyle%
\int_{E\cap \left( \ell _{m+1},\ell _{m}\right) }\left\Vert u\left(
\cdot
,t\right) \right\Vert _{L^{1}\left( B_{r}\right) }dt\text{ .}%
\end{array}
\tag{2.2.3}  \label{2.2.3}
\end{equation}%
Next, we take $m=2m^{\prime }$ and then sum (\ref{2.2.3})\ from
$m^{\prime }=1$ to infinity to deduce that
\begin{equation}
\begin{array}{ll}
& \quad \sum_{m^{\prime }=1}^{\infty }\left[ e^{-\eta \left( \gamma
+2\right) z^{2m^{\prime }}}\left\Vert u\left( \cdot ,\ell
_{2m^{\prime }}\right) \right\Vert _{L^{2}\left( \Omega \right)
}-e^{-\eta \left( \gamma +2\right) z^{2m^{\prime }+2}}\left\Vert
u\left( \cdot ,\ell _{2m^{\prime
}+2}\right) \right\Vert _{L^{2}\left( \Omega \right) }\right] \\
& \leq e^{dK\left( T,\left\Vert a\right\Vert ,\left\Vert
b\right\Vert \right) \beta \left( r,T,\left\Vert b\right\Vert
\right) }\sum_{m^{\prime }=1}^{\infty }\displaystyle\int_{E\cap
\left( \ell _{2m^{\prime }+1},\ell _{2m^{\prime }}\right)
}\left\Vert u\left( \cdot ,t\right) \right\Vert
_{L^{1}\left( B_{r}\right) }dt \\
& \leq e^{dK\left( T,\left\Vert a\right\Vert ,\left\Vert
b\right\Vert
\right) \beta \left( r,T,\left\Vert b\right\Vert \right) }\displaystyle%
\int_{E}\left\Vert u\left( \cdot ,t\right) \right\Vert _{L^{1}\left(
B_{r}\right) }dt\text{ .}%
\end{array}
\tag{2.2.4}  \label{2.2.4}
\end{equation}%
Since $e^{-\eta \left( \gamma +2\right) z^{2m^{\prime }+2}}$ tends
to zero as $m^{\prime }\rightarrow +\infty $, it holds that
\begin{equation}
\begin{array}{ll}
& \quad \sum_{m^{\prime }=1}^{\infty }\left[ e^{-\eta \left( \gamma
+2\right) z^{2m^{\prime }}}\left\Vert u\left( \cdot ,\ell
_{2m^{\prime }}\right) \right\Vert _{L^{2}\left( \Omega \right)
}-e^{-\eta \left( \gamma +2\right) z^{2m^{\prime }+2}}\left\Vert
u\left( \cdot ,\ell _{2m^{\prime }+2}\right) \right\Vert
_{L^{2}\left( \Omega \right) }\right] _{L^{2}\left(
\Omega \right) } \\
& =e^{-\eta \left( \gamma +2\right) z^{2}}\left\Vert u\left( \cdot
,\ell
_{2}\right) \right\Vert _{L^{2}\left( \Omega \right) }\text{ .}%
\end{array}
\tag{2.2.5}  \label{2.2.5}
\end{equation}%
Besides, one can easily check that
\begin{equation}
\begin{array}{ll}
\eta \left( \gamma +2\right) z^{2} & =\left( 1+d\beta \left(
r,T,\left\Vert b\right\Vert \right) \right) \left[ \frac{1}{\ell
_{1}-\ell }\right] \left(
\gamma +2\right) \sqrt{\gamma +2}\left( \sqrt{\gamma +2}+\sqrt{\gamma +1}%
\right) \\
& =C_{\left( \Omega ,n,q\right) }\frac{1}{\ell _{1}-\ell }\left[
\beta
\left( r,T,\left\Vert b\right\Vert \right) \right] ^{3}\text{ .}%
\end{array}
\tag{2.2.6}  \label{2.2.6}
\end{equation}%
Now, it follows from (\ref{2.2.4}), (\ref{2.2.5}) and (\ref{2.2.6})
that
\begin{equation*}
\left\Vert u\left( \cdot ,\ell _{2}\right) \right\Vert _{L^{2}\left(
\Omega
\right) }\leq e^{C_{\left( \Omega ,n,q\right) }\frac{1}{\ell _{1}-\ell }%
\left[ \beta \left( r,T,\left\Vert b\right\Vert \right) \right]
^{3}}e^{dK\left( T,\left\Vert a\right\Vert ,\left\Vert b\right\Vert
\right) \beta \left( r,T,\left\Vert b\right\Vert \right)
}\int_{E}\left\Vert u\left( \cdot ,t\right) \right\Vert
_{L^{1}\left( B_{r}\right) }dt\text{ .}
\end{equation*}%
This, along with the fact that
\begin{equation*}
\left\Vert u\left( \cdot ,T\right) \right\Vert _{L^{2}\left( \Omega
\right) }\leq e^{C_{0}T\left[ \left\Vert a\right\Vert _{L^{\infty
}\left( 0,T,L^{q}\left( \Omega \right) \right) }^{2}+\left\Vert
b\right\Vert _{L^{\infty }\left( \Omega \times \left( 0,T\right)
\right) }^{2}\right] }\left\Vert u\left( \cdot ,\ell _{2}\right)
\right\Vert _{L^{2}\left( \Omega \right) }\text{ ,}
\end{equation*}%
indicates that
\begin{equation*}
\left\Vert u\left( \cdot ,T\right) \right\Vert _{L^{2}\left( \Omega
\right) }\leq e^{\left( C_{0}+d\beta \left( r,T,\left\Vert
b\right\Vert \right) \right) K\left( T,\left\Vert a\right\Vert
,\left\Vert b\right\Vert \right) }e^{C_{\left( \Omega ,n,q\right)
}\frac{1}{\ell _{1}-\ell }\left[ \beta \left( r,T,\left\Vert
b\right\Vert \right) \right] ^{3}}\int_{E}\left\Vert u\left( \cdot
,t\right) \right\Vert _{L^{1}\left( B_{r}\right) }dt\text{ .}
\end{equation*}%
This leads to the desired results and completes the proof of Theorem
\ref{theorem1.1}.

\subsection{Proof of Proposition \ref{proposition2.2}}

We begin with introducing two quantities $G_{\lambda }$ and
$N_{\lambda ,\varphi }$ as follows. Let $x_{0}$ be the center of
$B_{r}$. Let $L\in \left( 0,T\right] $. For each $\lambda >0$, we
define
\begin{equation*}
G_{\lambda }\left( x,t\right) =\frac{1}{\left( L-t+\lambda \right) ^{n/2}}%
e^{-\frac{\left\vert x-x_{0}\right\vert ^{2}}{4\left( L-t+\lambda \right) }}%
\text{ ,\quad }\left( x,t\right) \in \mathbb{R}^{n}\times \left[
0,L\right] \text{ .}
\end{equation*}%
It is clear that $G_{\lambda }$ is a smooth function and satisfies
\begin{equation}
\left( \partial _{t}+\Delta \right) G_{\lambda }\left( x,t\right)
=0\text{
,\quad }\left( x,t\right) \in \mathbb{R}^{n}\times \left[ 0,L\right] \text{ .%
}  \tag{2.3.1}  \label{2.3.1}
\end{equation}%
Moreover, it holds that
\begin{equation}
\begin{array}{ll}
& \quad \frac{1}{2}\frac{d}{dt}\displaystyle\int_{\Omega }\left\vert
u\left(
x,t\right) \right\vert ^{2}G_{\lambda }\left( x,t\right) dx+\displaystyle%
\int_{\Omega }\left\vert \nabla u\left( x,t\right) \right\vert
^{2}G_{\lambda }\left( x,t\right) dx \\
& =\displaystyle\int_{\Omega }u\left( x,t\right) \left( \partial
_{t}-\Delta
\right) u\left( x,t\right) G_{\lambda }\left( x,t\right) dx\text{ ,}%
\end{array}
\tag{2.3.2}  \label{2.3.2}
\end{equation}%
for any $t\in \left( 0,L\right] $. This can be proved by a direct
computation. Also it can be derived from the following observation.
The quantity
\begin{equation*}
\int_{\Omega }\left( \partial _{t}-\Delta \right) \left( \left\vert
u\left( x,t\right) \right\vert ^{2}\right) G\left( x,t\right)
dx+\int_{\Omega }\left\vert u\left( x,t\right) \right\vert
^{2}\left( \partial _{t}+\Delta \right) G\left( x,t\right) dx
\end{equation*}%
where $G\in C^{\infty }$, has two expressions%
\begin{equation*}
\int_{\Omega }\frac{d}{dt}\left( \left\vert u\left( x,t\right)
\right\vert ^{2}G\left( x,t\right) \right) dx-\int_{\partial \Omega
}\left[ \partial _{\nu }\left( \left\vert u\left( x,t\right)
\right\vert ^{2}\right) G\left( x,t\right) -\left\vert u\left(
x,t\right) \right\vert ^{2}\partial _{\nu }G\left( x,t\right)
\right] d\sigma
\end{equation*}%
and
\begin{equation*}
\begin{array}{ll}
& \quad 2\displaystyle\int_{\Omega }\left[ u\left( x,t\right) \left(
\partial _{t}-\Delta \right) u\left( x,t\right) -\left\vert \nabla u\left(
x,t\right) \right\vert ^{2}\right] G\left( x,t\right) dx \\
& +\displaystyle\int_{\Omega }\left\vert u\left( x,t\right)
\right\vert
^{2}\left( \partial _{t}+\Delta \right) G\left( x,t\right) dx\text{ .}%
\end{array}%
\end{equation*}%
Because of (\ref{2.3.1}) and since $u=0$ on $\partial \Omega $, (\ref{2.3.2}%
) follows from the above two expressions with $G=G_{\lambda }$.

Next, we define, for each $\lambda >0$ and each $\varphi $ such that $%
\varphi \in C\left( \left[ \tau ,L\right] ;H^{1}\left( \Omega
\right) \right) $ for any $\tau \in \left( 0,L\right) $,
\begin{equation*}
N_{\lambda ,\varphi }\left( t\right) =\frac{\int_{\Omega }\left\vert
\nabla
\varphi \left( x,t\right) \right\vert ^{2}G_{\lambda }\left( x,t\right) dx}{%
\int_{\Omega }\left\vert \varphi \left( x,t\right) \right\vert
^{2}G_{\lambda }\left( x,t\right) dx}\text{ ,}
\end{equation*}%
where $t$ is in the set $\left\{ t\in \left( 0,L\right] ;\varphi
\left( \cdot ,t\right) \neq 0\text{ in }L^{2}\left( \Omega \right)
\right\} $.

\textbf{Proof of (\ref{2.1.3}) in Proposition \ref{proposition2.2}.}
\textit{The first step to prove (\ref{2.1.3}) is to estimate }$\frac{d}{dt}%
N_{\lambda ,u}\left( t\right) $\textit{. The desired estimate is a
consequence of the following lemma.}

\begin{lem} \label{lemma2.3} Let $\left( \varphi _{0},g\right) \in L^{2}\left( \Omega
\right) \times L^{2}\left( \Omega \times \left( 0,L\right) \right)
$\ and $\varphi =\varphi \left( x,t\right) $\ be the solution of
\begin{equation*}
\left\{
\begin{array}{ll}
\partial _{t}\varphi -\Delta \varphi =g & \quad \text{in~}\Omega
\times \left( 0,L\right) \ \text{,} \\
\varphi =0 & \quad \text{on~}\partial \Omega \times \left(
0,L\right) \text{ ,} \\
\varphi \left( \cdot ,0\right) =\varphi _{0}\text{ .} &
\end{array}%
\right.
\end{equation*}%
Then on the set $\left\{ t\in \left( 0,L\right] ;\varphi \left(
\cdot ,t\right) \neq 0\text{ in }L^{2}\left( \Omega \right) \right\}
$, the function $t\mapsto N_{\lambda ,\varphi }\left( t\right) $\ is
differentiable. Furthermore, it holds that
\begin{equation}
\frac{d}{dt}N_{\lambda ,\varphi }\left( t\right) \leq \frac{1}{L-t+\lambda }%
N_{\lambda ,\varphi }\left( t\right) +\frac{\int_{\Omega }\left\vert
g\left( x,t\right) \right\vert ^{2}G_{\lambda }\left( x,t\right)
dx}{\int_{\Omega }\left\vert \varphi \left( x,t\right) \right\vert
^{2}G_{\lambda }\left( x,t\right) dx}\text{ .}  \tag{2.3.3}
\label{2.3.3}
\end{equation}
\end{lem}

Lemma \ref{lemma2.3} is a direct consequence of estimate (3.26) in
\cite{PW}. We omit the proof.

\textit{The second step to prove (\ref{2.1.3}) is to estimate
}$\lambda
N_{\lambda ,u}\left( L\right) $\textit{\ by making use of (\ref{2.3.3}) and (%
\ref{2.3.2}). The desired estimate is stated in the following
lemma.}

\begin{lem} \label{lemma2.4} There exists a $C_{\left( \Omega ,n,q\right) }$\ such
that any non-trivial solution $u$\ to (\ref{1.1}) satisfies
\begin{equation*}
\begin{array}{ll}
& \quad \lambda N_{\lambda ,u}\left( L\right) +\frac{n}{4} \\
& \leq 8\left( \frac{\lambda }{L}+n\right) e^{2L\left( 1+\left\Vert
b\right\Vert _{L^{\infty }\left( \Omega \times \left( 0,L\right)
\right)
}^{2}\right) } \\
& \quad \times \log \left[ e^{\left( 1+\left( C_{\left( \Omega
,n,q\right) }+C_{0}\right) \left[ A\left( L,\left\Vert a\right\Vert
\right) +L\left\Vert b\right\Vert _{L^{\infty }\left( \Omega \times
\left(
0,L\right) \right) }^{2}\right] +\frac{m_{0}}{2L}\right) }\frac{\displaystyle%
\int_{\Omega }\left\vert u\left( x,0\right) \right\vert ^{2}dx}{\displaystyle%
\int_{\Omega }\left\vert u\left( x,L\right) \right\vert
^{2}dx}\right] \text{
,}%
\end{array}%
\end{equation*}%
where $m_{0}=\underset{x\in \Omega }{\sup }\left\vert
x-x_{0}\right\vert ^{2}$\ and $C_{0}$\ is given in (\ref{1.2}).
\end{lem}

\textbf{Proof of Lemma \ref{lemma2.4}.} Clearly, the solution $u$ to
(\ref{1.1}) holds the property that $u\in L^{2}\left( \tau
,T;H^{2}\cap H_{0}^{1}\left( \Omega \right) \right) \cap C\left(
\left[ \tau ,T\right] ;H_{0}^{1}\left( \Omega \right) \right) $ and
$\partial _{t}u\in L^{2}\left( \tau ,T;L^{2}\left( \Omega \right)
\right) $ for any $\tau \in \left( 0,L\right) $. One can easily
check that $N_{\lambda ,u}\left( t\right) $ is well-defined for any
$t\in \left( 0,L\right] $. We carry out the rest of the proof by
three steps as follows.

\textit{Step 1} .- We claim that for any $t\in \left( 0,L\right] $,%
\begin{equation}
\begin{array}{ll}
& \quad \frac{\lambda }{L+\lambda }e^{-2L\left( 1+\left\Vert
b\right\Vert _{L^{\infty }\left( \Omega \times \left( 0,L\right)
\right) }^{2}\right)
}N_{\lambda ,u}\left( L\right) \\
& \leq N_{\lambda ,u}\left( t\right) +C_{\left( \Omega ,n,q\right)
}L\left( \left\Vert a\right\Vert _{L^{\infty }\left( 0,T;L^{q}\left(
\Omega \right)
\right) }\right) ^{\frac{4}{2-p}}+C_{\left( \Omega ,n,q\right) }\frac{1}{%
L^{p-1}}\left\Vert a\right\Vert _{L^{\infty }\left( 0,T;L^{q}\left(
\Omega
\right) \right) }^{2}\text{ .}%
\end{array}
\tag{2.3.4}  \label{2.3.4}
\end{equation}%
To this ends, we apply Lemma \ref{lemma2.3} to $\left( \varphi
_{0},g\right) =\left( u\left( \cdot ,0\right) ,-au-b\cdot \nabla
u\right) $ and use Cauchy-Schwarz inequality to get that
\begin{equation}
\begin{array}{ll}
\frac{d}{dt}N_{\lambda ,u}\left( t\right) & \leq \frac{1}{L-t+\lambda }%
N_{\lambda ,u}\left( t\right) \\
& \quad +2\frac{\displaystyle\int_{\Omega }\left\vert au\left(
x,t\right) \right\vert ^{2}G_{\lambda }\left( x,t\right)
dx}{\displaystyle\int_{\Omega }\left\vert u\left( x,t\right)
\right\vert ^{2}G_{\lambda }\left( x,t\right) dx}+2\left\Vert
b\right\Vert _{L^{\infty }\left( \Omega \times \left(
0,L\right) \right) }^{2}N_{\lambda ,u}\left( t\right) \text{ .}%
\end{array}
\tag{2.3.5}  \label{2.3.5}
\end{equation}%
Since (\ref{A.2.1}) in Appendix holds,
\begin{equation}
\begin{array}{ll}
& \quad \displaystyle\int_{\Omega }\left\vert a\left( x,t\right)
u\left(
x,t\right) \right\vert ^{2}G_{\lambda }\left( x,t\right) dx \\
& \leq N_{\lambda ,u}\left( t\right) \displaystyle\int_{\Omega
}\left\vert
u\left( x,t\right) \right\vert ^{2}G_{\lambda }\left( x,t\right) dx \\
& \quad +C_{\left( \Omega ,n,q\right) }\left( \left( \left\Vert
a\right\Vert
_{L^{\infty }\left( 0,L;L^{q}\left( \Omega \right) \right) }\right) ^{\frac{4%
}{2-p}}+\frac{\left\Vert a\right\Vert _{L^{\infty }\left(
0,L;L^{q}\left(
\Omega \right) \right) }^{2}}{\left( L-t+\lambda \right) ^{p}}\right) %
\displaystyle\int_{\Omega }\left\vert u\left( x,t\right) \right\vert
^{2}G_{\lambda }\left( x,t\right) dx\text{ .}%
\end{array}
\tag{2.3.6}  \label{2.3.6}
\end{equation}%
It follows from (\ref{2.3.5}) and (\ref{2.3.6}) that
\begin{equation*}
\begin{array}{ll}
& \quad \frac{d}{dt}\left[ \left( L-t+\lambda \right) e^{-2t\left(
1+\left\Vert b\right\Vert _{L^{\infty }\left( \Omega \times \left(
0,L\right) \right) }^{2}\right) }N_{\lambda ,u}\left( t\right) \right] \\
& \leq C_{\left( \Omega ,n,q\right) }\left( \left\Vert a\right\Vert
_{L^{\infty }\left( 0,L;L^{q}\left( \Omega \right) \right) }\right) ^{\frac{4%
}{2-p}}\left( L-t+\lambda \right) e^{-2t\left( 1+\left\Vert
b\right\Vert _{L^{\infty }\left( \Omega \times \left( 0,L\right)
\right) }^{2}\right) }
\\
& \quad +C_{\left( \Omega ,n,q\right) }\left\Vert a\right\Vert
_{L^{\infty }\left( 0,L;L^{q}\left( \Omega \right) \right)
}^{2}\frac{1}{\left( L-t+\lambda \right) ^{p-1}}e^{-2t\left(
1+\left\Vert b\right\Vert
_{L^{\infty }\left( \Omega \times \left( 0,L\right) \right) }^{2}\right) }%
\text{ .}%
\end{array}%
\end{equation*}%
Integrating it over $\left[ t,L\right] $ with $t\in \left(
0,L\right) $, after some simple computations, we get (\ref{2.3.4}).

\textit{Step 2 }.- We claim that for any $t\in \left( 0,L/2\right]
$,
\begin{equation}
\begin{array}{ll}
& \quad \frac{d}{dt}\displaystyle\int_{\Omega }\left\vert u\left(
x,t\right) \right\vert ^{2}G_{\lambda }\left( x,t\right)
dx+\frac{1}{2}N_{\lambda ,u}\left( t\right)
\displaystyle\int_{\Omega }\left\vert u\left( x,t\right)
\right\vert ^{2}G_{\lambda }\left( x,t\right) dx \\
& \leq \left( C_{\left( \Omega ,n,q\right) }\left\Vert a\right\Vert
_{L^{\infty }\left( 0,L;L^{q}\left( \Omega \right) \right)
}^{2}+\left\Vert b\right\Vert _{L^{\infty }\left( \Omega \times
\left( 0,L\right) \right) }^{2}\right) \displaystyle\int_{\Omega
}\left\vert u\left( x,t\right)
\right\vert ^{2}G_{\lambda }\left( x,t\right) dx \\
& \quad +C_{\left( \Omega ,n,q\right) }\frac{1}{L}\left\Vert
a\right\Vert
_{L^{\infty }\left( 0,L;L^{q}\left( \Omega \right) \right) }\displaystyle%
\int_{\Omega }\left\vert u\left( x,t\right) \right\vert
^{2}G_{\lambda
}\left( x,t\right) dx\text{ .}%
\end{array}
\tag{2.3.7}  \label{2.3.7}
\end{equation}%
For this purpose, we first observe that (\ref{2.3.2}) is equivalent
to the following equality:
\begin{equation*}
\begin{array}{ll}
& \quad \frac{d}{dt}\displaystyle\int_{\Omega }\left\vert u\left(
x,t\right) \right\vert ^{2}G_{\lambda }\left( x,t\right)
dx+2N_{\lambda ,u}\left( t\right) \displaystyle\int_{\Omega
}\left\vert u\left( x,t\right)
\right\vert ^{2}G_{\lambda }\left( x,t\right) dx \\
& =2\displaystyle\int_{\Omega }u\left( x,t\right) \left( \partial
_{t}-\Delta \right) u\left( x,t\right) G_{\lambda }\left( x,t\right)
dx\text{
,}%
\end{array}%
\end{equation*}%
for any $t\in \left( 0,L\right] $. By this and by Cauchy-Schwarz
inequality,
it follows that%
\begin{equation}
\begin{array}{ll}
& \quad \frac{d}{dt}\displaystyle\int_{\Omega }\left\vert u\left(
x,t\right) \right\vert ^{2}G_{\lambda }\left( x,t\right)
dx+N_{\lambda ,u}\left( t\right) \displaystyle\int_{\Omega
}\left\vert u\left( x,t\right)
\right\vert ^{2}G_{\lambda }\left( x,t\right) dx \\
& \leq 2\displaystyle\int_{\Omega }\left\vert a\left( x,t\right)
\right\vert \left\vert u\left( x,t\right) \right\vert ^{2}G_{\lambda
}\left( x,t\right)
dx \\
& \quad +\left\Vert b\right\Vert _{L^{\infty }\left( \Omega \times
\left( 0,L\right) \right) }^{2}\displaystyle\int_{\Omega }\left\vert
u\left(
x,t\right) \right\vert ^{2}G_{\lambda }\left( x,t\right) dx\text{ .}%
\end{array}
\tag{2.3.8}  \label{2.3.8}
\end{equation}%
Because of (\ref{A.2.2}) in Appendix,
\begin{equation}
\begin{array}{ll}
& \quad \displaystyle\int_{\Omega }a\left( x,t\right) \left\vert
u\left(
x,t\right) \right\vert ^{2}G_{\lambda }\left( x,t\right) dx \\
& \leq \frac{1}{4}\displaystyle\int_{\Omega }\left\vert \nabla
u\left(
x,t\right) \right\vert ^{2}G_{\lambda }\left( x,t\right) dx \\
& \quad +C_{\left( \Omega ,n,q\right) }\left\Vert a\right\Vert
_{L^{\infty }\left( 0,L;L^{q}\left( \Omega \right) \right) }\left(
\left\Vert
a\right\Vert _{L^{\infty }\left( 0,L;L^{q}\left( \Omega \right) \right) }+%
\frac{1}{L-t+\lambda }\right) \displaystyle\int_{\Omega }\left\vert
u\left(
x,t\right) \right\vert ^{2}G_{\lambda }\left( x,t\right) dx\text{ .}%
\end{array}
\tag{2.3.9}  \label{2.3.9}
\end{equation}%
We directly get (\ref{2.3.7}) from (\ref{2.3.8}) and (\ref{2.3.9}).

\textit{Step 3 }.- Conclusion. By (\ref{2.3.4}) and (\ref{2.3.7}),
we deduce that for any $t\in \left( 0,L/2\right] $,
\begin{equation}
\begin{array}{ll}
& \quad \frac{d}{dt}\displaystyle\int_{\Omega }\left\vert u\left(
x,t\right)
\right\vert ^{2}G_{\lambda }\left( x,t\right) dx \\
& \quad +\frac{1}{2}\frac{\lambda }{L+\lambda }e^{-2L\left(
1+\left\Vert b\right\Vert _{L^{\infty }\left( \Omega \times \left(
0,L\right) \right) }^{2}\right) }N_{\lambda ,u}\left( L\right)
\displaystyle\int_{\Omega }\left\vert u\left( x,t\right) \right\vert
^{2}G_{\lambda }\left( x,t\right)
dx \\
& \leq \left( C_{\left( \Omega ,n,q\right) }\left\Vert a\right\Vert
_{L^{\infty }\left( 0,L;L^{q}\left( \Omega \right) \right)
}^{2}+\left\Vert b\right\Vert _{L^{\infty }\left( \Omega \times
\left( 0,L\right) \right) }^{2}\right) \displaystyle\int_{\Omega
}\left\vert u\left( x,t\right)
\right\vert ^{2}G_{\lambda }\left( x,t\right) dx \\
& \quad +C_{\left( \Omega ,n,q\right) }\frac{1}{L}\left\Vert
a\right\Vert
_{L^{\infty }\left( 0,L;L^{q}\left( \Omega \right) \right) }\displaystyle%
\int_{\Omega }\left\vert u\left( x,t\right) \right\vert
^{2}G_{\lambda
}\left( x,t\right) dx\text{ } \\
& \quad +C_{\left( \Omega ,n,q\right) }\left[ L\left( \left\Vert
a\right\Vert _{L^{\infty }\left( 0,L;L^{q}\left( \Omega \right)
\right) }\right) ^{\frac{4}{2-p}}+\frac{1}{L^{p-1}}\left\Vert
a\right\Vert
_{L^{\infty }\left( 0,L;L^{q}\left( \Omega \right) \right) }^{2}\right] \\
& \quad \times \displaystyle\int_{\Omega }\left\vert u\left(
x,t\right)
\right\vert ^{2}G_{\lambda }\left( x,t\right) dx\text{ .}%
\end{array}
\tag{2.3.10}  \label{2.3.10}
\end{equation}%
Recall that
\begin{equation*}
\begin{array}{ll}
\frac{1}{L}A\left( L,\left\Vert a\right\Vert \right) & =\left\Vert
a\right\Vert _{L^{\infty }\left( 0,L;L^{q}\left( \Omega \right)
\right) }^{2}+\frac{1}{L}\left\Vert a\right\Vert _{L^{\infty }\left(
0,L;L^{q}\left(
\Omega \right) \right) } \\
& \quad +L\left( \left\Vert a\right\Vert _{L^{\infty }\left(
0,L;L^{q}\left( \Omega \right) \right) }\right)
^{\frac{4}{2-p}}+\frac{1}{L^{p-1}}\left\Vert
a\right\Vert _{L^{\infty }\left( 0,L;L^{q}\left( \Omega \right) \right) }^{2}%
\text{ .}%
\end{array}%
\end{equation*}%
This, together with (\ref{2.3.10}), gives that%
\begin{equation*}
\begin{array}{ll}
& \quad \frac{d}{dt}\displaystyle\int_{\Omega }\left\vert u\left(
x,t\right)
\right\vert ^{2}G_{\lambda }\left( x,t\right) dx \\
& \quad +\frac{1}{2}\frac{\lambda }{L+\lambda }e^{-2L\left(
1+\left\Vert b\right\Vert _{L^{\infty }\left( \Omega \times \left(
0,L\right) \right) }^{2}\right) }N_{\lambda ,u}\left( L\right)
\displaystyle\int_{\Omega }\left\vert u\left( x,t\right) \right\vert
^{2}G_{\lambda }\left( x,t\right)
dx \\
& \leq \left( C_{\left( \Omega ,n,q\right) }\frac{1}{L}A\left(
L,\left\Vert a\right\Vert \right) +\left\Vert b\right\Vert
_{L^{\infty }\left( \Omega \times \left( 0,L\right) \right)
}^{2}\right) \displaystyle\int_{\Omega }\left\vert u\left(
x,t\right) \right\vert ^{2}G_{\lambda }\left( x,t\right)
dx\text{ .}%
\end{array}%
\end{equation*}%
From this, we deduce that for any $t\in \left( 0,L/2\right] $,%
\begin{equation*}
\begin{array}{ll}
\displaystyle\int_{\Omega }\left\vert u\left( x,t\right) \right\vert
^{2}G_{\lambda }\left( x,t\right) dx & \leq \text{exp}\left( -t\left[ \frac{1%
}{2}\frac{\lambda }{L+\lambda }e^{-2L\left( 1+\left\Vert
b\right\Vert _{L^{\infty }\left( \Omega \times \left( 0,L\right)
\right) }^{2}\right)
}N_{\lambda }\left( u,L\right) \right] \right) \\
& \quad \times \text{exp}\left( t\left[ C_{\left( \Omega ,n,q\right) }\frac{1%
}{L}A\left( L,\left\Vert a\right\Vert \right) +\left\Vert
b\right\Vert _{L^{\infty }\left( \Omega \times \left( 0,L\right)
\right) }^{2}\right]
\right) \\
& \quad \times \displaystyle\int_{\Omega }\left\vert u\left(
x,0\right)
\right\vert ^{2}G_{\lambda }\left( x,0\right) dx\text{ .}%
\end{array}%
\end{equation*}%
Taking $t=L/2$ in the above, we see that%
\begin{equation}
\begin{array}{ll}
& \quad \frac{1}{\left( L/2+\lambda \right)
^{n/2}}\displaystyle\int_{\Omega }\left\vert u\left( x,L/2\right)
\right\vert ^{2}e^{-\frac{\left\vert
x-x_{0}\right\vert ^{2}}{4\left( L/2+\lambda \right) }}dx \\
& \leq \text{exp}\left( -\left[ \frac{\lambda L}{4\left( L+\lambda \right) }%
e^{-2L\left( 1+\left\Vert b\right\Vert _{L^{\infty }\left( \Omega
\times \left( 0,L\right) \right) }^{2}\right) }N_{\lambda }\left(
u,L\right) \right]
\right) \\
& \quad \times \text{exp}\left( C_{\left( \Omega ,n,q\right)
}A\left( L,\left\Vert a\right\Vert \right) +L\left\Vert b\right\Vert
_{L^{\infty
}\left( \Omega \times \left( 0,L\right) \right) }^{2}\right) \\
& \quad \times \frac{1}{\left( L+\lambda \right) ^{n/2}}\displaystyle%
\int_{\Omega }\left\vert u\left( x,0\right) \right\vert ^{2}e^{-\frac{%
\left\vert x-x_{0}\right\vert ^{2}}{4\left( L+\lambda \right) }}dx\text{ .}%
\end{array}
\tag{2.3.11}  \label{2.3.11}
\end{equation}%
On the other hand, it is clear that
\begin{equation*}
\begin{array}{ll}
& \quad \displaystyle\int_{\Omega }\left\vert u\left( x,L\right)
\right\vert
^{2}dx \\
& \leq e^{C_{0}L\left( \left\Vert a\right\Vert _{L^{\infty }\left(
0,L;L^{q}\left( \Omega \right) \right) }^{2}+\left\Vert b\right\Vert
_{L^{\infty }\left( \Omega \times \left( 0,L\right) \right) }^{2}\right) }%
\displaystyle\int_{\Omega }\left\vert u\left( x,L/2\right)
\right\vert ^{2}dx
\\
& \leq e^{C_{0}L\left( \left\Vert a\right\Vert _{L^{\infty }\left(
0,L;L^{q}\left( \Omega \right) \right) }^{2}+\left\Vert b\right\Vert
_{L^{\infty }\left( \Omega \times \left( 0,L\right) \right) }^{2}\right) }e^{%
\frac{m_{0}}{2L}}\displaystyle\int_{\Omega }\left\vert u\left(
x,L/2\right) \right\vert ^{2}e^{-\frac{\left\vert x-x_{0}\right\vert
^{2}}{4\left(
L/2+\lambda \right) }}dx\text{ .}%
\end{array}%
\end{equation*}%
This, together with (\ref{2.3.11}), yields that%
\begin{equation*}
\begin{array}{ll}
\displaystyle\int_{\Omega }\left\vert u\left( x,L\right) \right\vert
^{2}dx
& \leq \text{exp}\left( -\left[ \frac{\lambda L}{4\left( L+\lambda \right) }%
e^{-2L\left( 1+\left\Vert b\right\Vert _{L^{\infty }\left( \Omega
\times \left( 0,L\right) \right) }^{2}\right) }N_{\lambda }\left(
u,L\right) \right]
\right) \\
& \quad \times \text{exp}\left( \left( C_{\left( \Omega ,n,q\right)
}+C_{0}\right) \left[ A\left( L,\left\Vert a\right\Vert \right)
+L\left\Vert b\right\Vert _{L^{\infty }\left( \Omega \times \left(
0,L\right) \right)
}^{2}\right] +\frac{m_{o}}{2L}\right) \\
& \quad \times \displaystyle\int_{\Omega }\left\vert u\left(
x,0\right)
\right\vert ^{2}dx\text{ ,}%
\end{array}%
\end{equation*}%
from which it follows that%
\begin{equation}
\begin{array}{ll}
\lambda N_{\lambda ,u}\left( L\right) & \leq 4\left( \frac{\lambda }{L}%
+1\right) e^{2L\left( 1+\left\Vert b\right\Vert _{L^{\infty }\left(
\Omega
\times \left( 0,L\right) \right) }^{2}\right) } \\
& \quad \times \text{log}\left[ e^{\left( \left( C_{\left( \Omega
,n,q\right) }+C_{0}\right) \left[ A\left( L,\left\Vert a\right\Vert
\right) +L\left\Vert b\right\Vert _{L^{\infty }\left( \Omega \times
\left(
0,L\right) \right) }^{2}\right] +\frac{m_{0}}{2L}\right) }\frac{\displaystyle%
\int_{\Omega }\left\vert u\left( x,0\right) \right\vert ^{2}dx}{\displaystyle%
\int_{\Omega }\left\vert u\left( x,L\right) \right\vert
^{2}dx}\right] \text{
.}%
\end{array}
\tag{2.3.12}  \label{2.3.12}
\end{equation}%
Clearly, it holds that
\begin{equation}
\frac{n}{4}\leq \frac{n}{4}\text{log}\left[ e^{\left( 1+\left(
C_{\left( \Omega ,n,q\right) }+C_{0}\right) \left[ A\left(
L,\left\Vert a\right\Vert \right) +L\left\Vert b\right\Vert
_{L^{\infty }\left( \Omega \times \left( 0,L\right) \right)
}^{2}\right] +\frac{m}{2L}\right) }\frac{\int_{\Omega }\left\vert
u\left( x,0\right) \right\vert ^{2}dx}{\int_{\Omega }\left\vert
u\left( x,L\right) \right\vert ^{2}dx}\right] \text{ .}
\tag{2.3.13} \label{2.3.13}
\end{equation}%
Now, the desired estimate in Lemma \ref{lemma2.4} follows
immediately from (\ref{2.3.12}) and (\ref{2.3.13}). This completes
the proof of Lemma \ref{lemma2.4}.

\textit{The third step to prove (\ref{2.1.3}) is to get an estimate
of }$\int_{\Omega }\left\vert u\left( x,L\right) \right\vert
^{2}e^{-\frac{\left\vert x-x_{0}\right\vert ^{2}}{4\lambda
}}dx$\textit{\ in
term of }$\int_{B_{r}}\left\vert u\left( x,L\right) \right\vert ^{2}e^{-%
\frac{\left\vert x-x_{0}\right\vert ^{2}}{4\lambda }}dx$\textit{. It
is a consequence of the following lemma.}

\begin{lem} \label{lemma2.5} For any non-trivial $f\in H^{1}\left( \Omega \right) $\ and any $%
\lambda >0$, it holds that
\begin{equation*}
\begin{array}{ll}
\displaystyle\int_{\Omega }\left\vert f\left( x\right) \right\vert ^{2}e^{-%
\frac{\left\vert x-x_{0}\right\vert ^{2}}{4\lambda }}dx & \leq \displaystyle%
\int_{B_{r}}\left\vert f\left( x\right) \right\vert
^{2}e^{-\frac{\left\vert
x-x_{0}\right\vert ^{2}}{4\lambda }}dx \\
& \quad +\frac{16\lambda }{r^{2}}\left( \lambda N_{\lambda ,f}\left(
L\right) +\frac{n}{4}\right) \displaystyle\int_{\Omega }\left\vert
f\left(
x\right) \right\vert ^{2}e^{-\frac{\left\vert x-x_{0}\right\vert ^{2}}{%
4\lambda }}dx\text{ .}%
\end{array}%
\end{equation*}
\end{lem}

\textbf{Proof of Lemma \ref{lemma2.5}.} We first observe that
\begin{equation}
\begin{array}{ll}
& \quad \displaystyle\int_{\Omega }\left\vert f\left( x\right)
\right\vert
^{2}e^{-\frac{\left\vert x-x_{0}\right\vert ^{2}}{4\lambda }}dx \\
& \leq \displaystyle\int_{B_{r}}\left\vert f\left( x\right)
\right\vert
^{2}e^{-\frac{\left\vert x-x_{0}\right\vert ^{2}}{4\lambda }}dx+\displaystyle%
\int_{\Omega \cap \left\{ \left\vert x-x_{o}\right\vert \geq
r\right\} }\left\vert f\left( x\right) \right\vert
^{2}e^{-\frac{\left\vert
x-x_{0}\right\vert ^{2}}{4\lambda }}dx \\
& \leq \displaystyle\int_{B_{r}}\left\vert f\left( x\right)
\right\vert
^{2}e^{-\frac{\left\vert x-x_{0}\right\vert ^{2}}{4\lambda }}dx+\frac{%
16\lambda }{r^{2}}\displaystyle\int_{\Omega }\frac{\left\vert
x-x_{0}\right\vert ^{2}}{16\lambda }\left\vert f\left( x\right)
\right\vert
^{2}e^{-\frac{\left\vert x-x_{0}\right\vert ^{2}}{4\lambda }}dx\text{ .}%
\end{array}
\tag{2.3.14}  \label{2.3.14}
\end{equation}%
Next, we claim that%
\begin{equation}
\begin{array}{ll}
& \quad \displaystyle\int_{\Omega }\frac{\left\vert x-x_{0}\right\vert ^{2}}{%
16\lambda }\left\vert f\left( x\right) \right\vert
^{2}e^{-\frac{\left\vert
x-x_{0}\right\vert ^{2}}{4\lambda }}dx \\
& \leq \lambda \displaystyle\int_{\Omega }\left\vert \nabla f\left(
x\right)
\right\vert ^{2}e^{-\frac{\left\vert x-x_{0}\right\vert ^{2}}{4\lambda }}dx+%
\frac{n}{4}\displaystyle\int_{\Omega }\left\vert f\left( x\right)
\right\vert ^{2}e^{-\frac{\left\vert x-x_{0}\right\vert ^{2}}{4\lambda }}dx%
\text{ .}%
\end{array}
\tag{2.3.15}  \label{2.3.15}
\end{equation}%
When this is done, the desired estimate in Lemma \ref{lemma2.5} follows at once from (\ref%
{2.3.14}) and (\ref{2.3.15}). It remains to show (\ref{2.3.15}).
This can be done by what follows (see also \cite[page 211]{EFV}).
\begin{equation}
\begin{array}{ll}
& \quad \displaystyle\int_{\Omega }\left\vert x-x_{0}\right\vert
^{2}\left\vert f\left( x\right) \right\vert ^{2}e^{-\frac{\left\vert
x-x_{0}\right\vert ^{2}}{4\lambda }}dx \\
& =\displaystyle\int_{\Omega }\left( x-x_{0}\right) \left\vert
f\left(
x\right) \right\vert ^{2}\cdot \left( -2\lambda \right) \nabla e^{-\frac{%
\left\vert x-x_{0}\right\vert ^{2}}{4\lambda }}dx \\
& =-2\lambda \displaystyle\int_{\partial \Omega }\left( \left(
x-x_{0}\right) \cdot \nu \right) \left\vert f\left( x\right)
\right\vert ^{2}e^{-\frac{\left\vert x-x_{0}\right\vert
^{2}}{4\lambda }}d\sigma +2\lambda n\displaystyle\int_{\Omega
}\left\vert f\left( x\right) \right\vert ^{2}e^{-\frac{\left\vert
x-x_{0}\right\vert ^{2}}{4\lambda }}dx
\\
& \quad +4\lambda \displaystyle\int_{\Omega }\left( x-x_{0}\right)
f\left( x\right) \cdot \nabla f\left( x\right) e^{-\frac{\left\vert
x-x_{0}\right\vert ^{2}}{4\lambda }}dx \\
& \leq 2\lambda n\displaystyle\int_{\Omega }\left\vert f\left(
x\right) \right\vert ^{2}e^{-\frac{\left\vert x-x_{0}\right\vert
^{2}}{4\lambda }}dx
\\
& \quad +\frac{1}{2}\displaystyle\int_{\Omega }16\lambda
^{2}\left\vert \nabla f\left( x\right) \right\vert
^{2}e^{-\frac{\left\vert x-x_{0}\right\vert ^{2}}{4\lambda
}}dx+\frac{1}{2}\displaystyle\int_{\Omega }\left\vert
x-x_{0}\right\vert ^{2}\left\vert f\left( x\right) \right\vert
^{2}e^{-\frac{\left\vert x-x_{0}\right\vert ^{2}}{4\lambda }}dx\text{ .}%
\end{array}
\tag{2.3.16}  \label{2.3.16}
\end{equation}%
In (\ref{2.3.16}), it is used in the first equality that $\left(
-2\lambda
\right) \nabla e^{-\frac{\left\vert x-x_{0}\right\vert ^{2}}{4\lambda }%
}=\left( x-x_{0}\right) e^{-\frac{\left\vert x-x_{0}\right\vert ^{2}}{%
4\lambda }}$; integration by parts is applied in the second
equality; Cauchy-Schwarz inequality, along with the assumption that
$\Omega $ is convex, is applied in the last inequality. This
completes the proof of Lemma \ref{lemma2.5}.

\textit{The last step to prove (\ref{2.1.3}) of Proposition
\ref{proposition2.2} is to drop the
weight function }$e^{-\frac{\left\vert x-x_{0}\right\vert ^{2}}{4\lambda }}$%
\textit{\ in the integrands.}

Recall that for any $\tau \in \left( 0,L\right) $, $u\in L^{2}\left(
\tau ,T;H^{2}\cap H_{0}^{1}\left( \Omega \right) \right) \cap
C\left( \left[ \tau ,T\right] ;H_{0}^{1}\left( \Omega \right)
\right) $ and $\partial _{t}u\in L^{2}\left( \tau ,T;L^{2}\left(
\Omega \right) \right) $. Without a
loss of generality, we assume that $u$ is non-trivial in order that $%
N_{\lambda ,u}\left( t\right) $ is well-defined for any $t\in \left( 0,L%
\right] $. We apply Lemma \ref{lemma2.5} where $f=u\left( \cdot
,L\right) $ to get that
\begin{equation}
\begin{array}{ll}
\displaystyle\int_{\Omega }\left\vert u\left( x,L\right) \right\vert ^{2}e^{-%
\frac{\left\vert x-x_{0}\right\vert ^{2}}{4\lambda }}dx & \leq \displaystyle%
\int_{B_{r}}\left\vert u\left( x,L\right) \right\vert ^{2}e^{-\frac{%
\left\vert x-x_{0}\right\vert ^{2}}{4\lambda }}dx \\
& \quad +\frac{16\lambda }{r^{2}}\left( \lambda N_{\lambda ,u}\left(
L\right) +\frac{n}{4}\right) \displaystyle\int_{\Omega }\left\vert
u\left(
x,L\right) \right\vert ^{2}e^{-\frac{\left\vert x-x_{0}\right\vert ^{2}}{%
4\lambda }}dx\text{ .}%
\end{array}
\tag{2.3.17}  \label{2.3.17}
\end{equation}%
By Lemma \ref{lemma2.4}, it holds that
\begin{equation}
\lambda N_{\lambda ,u}\left( L\right) +\frac{n}{4}\leq
\frac{1}{16}\left( \frac{\lambda }{L}+n\right) Z_{u}\text{ ,}
\tag{2.3.18}  \label{2.3.18}
\end{equation}%
where%
\begin{equation}
\begin{array}{ll}
Z_{u} & =16\times 8e^{2L\left( 1+\left\Vert b\right\Vert _{L^{\infty
}\left(
\Omega \times \left( 0,L\right) \right) }^{2}\right) } \\
& \quad \times \text{log}\left[ e^{\left( 1+\left( C_{\left( \Omega
,n,q\right) }+C_{0}\right) \left[ A\left( L,\left\Vert a\right\Vert
\right) +L\left\Vert b\right\Vert _{L^{\infty }\left( \Omega \times
\left(
0,L\right) \right) }^{2}\right] +\frac{m_{_{0}}}{2L}\right) }\frac{%
\int_{\Omega }\left\vert u\left( x,0\right) \right\vert
^{2}dx}{\int_{\Omega
}\left\vert u\left( x,L\right) \right\vert ^{2}dx}\right] \text{ .}%
\end{array}
\tag{2.3.19}  \label{2.3.19}
\end{equation}%
Combining (\ref{2.3.17}) and (\ref{2.3.18}), we get that for any $\lambda >0$%
,
\begin{equation}
\begin{array}{ll}
\displaystyle\int_{\Omega }\left\vert u\left( x,L\right) \right\vert ^{2}e^{-%
\frac{\left\vert x-x_{o}\right\vert ^{2}}{4\lambda }}dx & \leq \displaystyle%
\int_{B_{r}}\left\vert u\left( x,L\right) \right\vert ^{2}e^{-\frac{%
\left\vert x-x_{0}\right\vert ^{2}}{4\lambda }}dx \\
& \quad +\frac{\lambda }{r^{2}}\left( \frac{\lambda }{L}+n\right) Z_{u}%
\displaystyle\int_{\Omega }\left\vert u\left( x,L\right) \right\vert ^{2}e^{-%
\frac{\left\vert x-x_{0}\right\vert ^{2}}{4\lambda }}dx\text{ .}%
\end{array}
\tag{2.3.20}  \label{2.3.20}
\end{equation}%
We take%
\begin{equation*}
\lambda =\frac{1}{2}\left( -nL+\sqrt{n^{2}L^{2}+\frac{2Lr^{2}}{Z_{u}}}%
\right) \text{ .}
\end{equation*}%
Clearly it solves
\begin{equation}
\frac{\lambda }{r^{2}}\left( \frac{\lambda }{L}+n\right) Z_{u}=\frac{1}{2}%
\text{ .}  \tag{2.3.21}  \label{2.3.21}
\end{equation}%
Then it follows from (\ref{2.3.20}) and (\ref{2.3.21}) that
\begin{equation}
\int_{\Omega }\left\vert u\left( x,L\right) \right\vert ^{2}dx\leq 2e^{\frac{%
m_{0}}{4\lambda }}\int_{B_{r}}\left\vert u\left( x,L\right)
\right\vert ^{2}dx\text{ ,}  \tag{2.3.22}  \label{2.3.22}
\end{equation}%
where $m_{0}$ is given by Lemma \ref{lemma2.4}.

On the other hand, it holds that
\begin{equation}
e^{\frac{m_{0}}{4\lambda }}\leq e^{\frac{\left( n+1\right) m_{0}}{2}\frac{%
Z_{u}}{r^{2}}}  \tag{2.3.23}  \label{2.3.23}
\end{equation}%
because
\begin{equation*}
\begin{array}{ll}
\frac{1}{\lambda }=\frac{2}{-nL+\sqrt{n^{2}L^{2}+\frac{2Lr^{2}}{Z_{u}}}} & =%
\frac{Z_{u}}{Lr^{2}}\left(
nL+\sqrt{n^{2}L^{2}+\frac{2Lr^{2}}{Z_{u}}}\right)
\\
& \leq \frac{Z_{u}}{Lr^{2}}\left( nL+\sqrt{n^{2}L^{2}+\frac{4L^{2}r^{2}}{%
m_{0}}}\right) \\
& \leq \frac{Z_{u}}{r^{2}}\left( n+\sqrt{n^{2}+4}\right) \text{ .}%
\end{array}%
\end{equation*}%
In the first inequality of the above, we used that
$Z_{u}>\frac{m_{0}}{2L}$. Now it follows from (\ref{2.3.22}) and
(\ref{2.3.23}) that
\begin{equation}
\int_{\Omega }\left\vert u\left( x,L\right) \right\vert ^{2}dx\leq 2e^{\frac{%
\left( n+1\right)
m_{0}}{2}\frac{Z_{u}}{r^{2}}}\int_{B_{r}}\left\vert u\left(
x,L\right) \right\vert ^{2}dx\text{ .}  \tag{2.3.24}  \label{2.3.24}
\end{equation}%
Next, by (\ref{2.3.19}), there is a $C=C_{(\Omega ,n,q)}>2$\ such that%
\begin{equation*}
\frac{\left( n+1\right) m_{0}}{2}\frac{Z_{u}}{r^{2}}\leq C\beta
\left( r,T,\left\Vert b\right\Vert \right) \text{log}\left[
e^{C\left( K\left(
T,\left\Vert a\right\Vert ,\left\Vert b\right\Vert \right) +\frac{1}{L}%
\right) }\frac{\int_{\Omega }\left\vert u\left( x,0\right) \right\vert ^{2}dx%
}{\int_{\Omega }\left\vert u\left( x,L\right) \right\vert
^{2}dx}\right] \text{ .}
\end{equation*}%
This, together with (\ref{2.3.24}), yields that
\begin{equation*}
\int_{\Omega }\left\vert u\left( x,L\right) \right\vert ^{2}dx\leq
2\left[ e^{C\left( K\left( T,\left\Vert a\right\Vert ,\left\Vert
b\right\Vert \right) +\frac{1}{L}\right) }\frac{\int_{\Omega
}\left\vert u\left( x,0\right) \right\vert ^{2}dx}{\int_{\Omega
}\left\vert u\left( x,L\right) \right\vert ^{2}dx}\right] ^{C\beta
\left( r,T,\left\Vert b\right\Vert
\right) }\int_{B_{r}}\left\vert u\left( x,L\right) \right\vert ^{2}dx\text{ .%
}
\end{equation*}%
In summary, we conclude that
\begin{equation*}
\begin{array}{ll}
\displaystyle\int_{\Omega }\left\vert u\left( x,L\right) \right\vert
^{2}dx & \leq \left( 2\displaystyle\int_{B_{r}}\left\vert u\left(
x,L\right) \right\vert ^{2}dx\right) ^{\frac{1}{1+C\beta \left(
r,T,\left\Vert
b\right\Vert \right) }} \\
& \quad \times \left( e^{C\left( K\left( T,\left\Vert a\right\Vert
,\left\Vert b\right\Vert \right) +\frac{1}{L}\right) }\displaystyle%
\int_{\Omega }\left\vert u\left( x,0\right) \right\vert ^{2}dx\right) ^{%
\frac{C\beta \left( r,T,\left\Vert b\right\Vert \right) }{1+C\beta
\left(
r,T,\left\Vert b\right\Vert \right) }}%
\end{array}%
\end{equation*}%
which leads to (\ref{2.1.3}).

\textbf{Proof of (\ref{2.1.4}) in Proposition \ref{proposition2.2}
.} Let $0\leq t_{1}<t_{2}\leq T$. The estimate (\ref{2.1.3}) implies that%
\begin{equation}
\begin{array}{ll}
\left\Vert u\left( \cdot ,t_{2}\right) \right\Vert _{L^{2}\left(
\Omega \right) } & \leq \left( \sqrt{C}\left\Vert u\left( \cdot
,t_{2}\right) \right\Vert _{L^{2}\left( B_{r/2}\right) }\right)
^{1-\alpha \left(
r/2,T,\left\Vert b\right\Vert \right) } \\
& \quad \times \left( e^{C\left( K\left( T,\left\Vert a\right\Vert
,\left\Vert b\right\Vert \right) +\frac{1}{t_{2}-t_{1}}\right)
}\left\Vert u\left( \cdot ,t_{1}\right) \right\Vert _{L^{2}\left(
\Omega \right)
}\right) ^{\alpha \left( r/2,T,\left\Vert b\right\Vert \right) }\text{ .}%
\end{array}
\tag{2.3.25}  \label{2.3.25}
\end{equation}%
On the other hand, by Nash inequality and Poincar\'{e} inequality,
there
exists $c>0$ (depending only on $\Omega $ and $n$) such that%
\begin{equation}
\left\Vert u\left( \cdot ,t_{2}\right) \right\Vert _{L^{2}\left(
B_{r/2}\right) }^{1+2/n}\leq \frac{c}{r}\left\Vert u\left( \cdot
,t_{2}\right) \right\Vert _{L^{1}\left( B_{r}\right)
}^{2/n}\left\Vert
\nabla u\left( \cdot ,t_{2}\right) \right\Vert _{L^{2}\left( \Omega \right) }%
\text{ .}  \tag{2.3.26}  \label{2.3.26}
\end{equation}%
It follows from the standard energy method that
\begin{equation}
\left\Vert \nabla u\left( \cdot ,t_{2}\right) \right\Vert
_{L^{2}\left(
\Omega \right) }\leq \frac{1}{\left( t_{2}-t_{1}\right) ^{1/2}}e^{c\left( 1+T%
\left[ \left\Vert a\right\Vert _{L^{\infty }\left( 0,T;L^{q}\left(
\Omega \right) \right) }^{2}+\left\Vert b\right\Vert _{L^{\infty
}\left( \Omega \times \left( 0,T\right) \right) }^{2}\right] \right)
}\left\Vert u\left( \cdot ,t_{1}\right) \right\Vert _{L^{2}\left(
\Omega \right) }\text{ ,} \tag{2.3.27}  \label{2.3.27}
\end{equation}%
where $c>0$ (depends only on $\Omega $, $n$ and $q$). Combining (\ref{2.3.25}%
), (\ref{2.3.26}) and (\ref{2.3.27}), we deduce that there is a
positive
constant $d$ (only depending on $\Omega $, $n$ and $q$) such that%
\begin{equation*}
\begin{array}{ll}
\left\Vert u\left( \cdot ,t_{2}\right) \right\Vert _{L^{2}\left(
\Omega \right) } & \leq \left( \frac{1}{r}\right) ^{\frac{1-\alpha
\left( r/2,T,\left\Vert b\right\Vert \right) }{1+2/n}}e^{d\left(
K\left(
T,\left\Vert a\right\Vert ,\left\Vert b\right\Vert \right) +\frac{1}{%
t_{2}-t_{1}}\right) } \\
& \quad \times \left\Vert u\left( \cdot ,t_{2}\right) \right\Vert
_{L^{1}\left( B_{r}\right) }^{\frac{2/n}{1+2/n}\left[ 1-\alpha
\left( r/2,T,\left\Vert b\right\Vert \right) \right] }\left\Vert
u\left( \cdot
,t_{1}\right) \right\Vert _{L^{2}\left( \Omega \right) }^{1-\frac{2/n}{1+2/n}%
\left[ 1-\alpha \left( r/2,T,\left\Vert b\right\Vert \right) \right]
}\text{
.}%
\end{array}%
\end{equation*}%
This, together with some simple computations, leads to estimate (\ref{2.1.4}%
), and completes the proof of Proposition \ref{proposition2.2}.

\section{Applications to bang-bang controls}

Throughout this section, we assume that $a\in L^{\infty }\left(
\Omega \times \left( 0,T\right) \right) $, $B\in L^{\infty }\left(
\Omega \times \left( 0,T\right) \right) ^{n}$ with div$B\in
L^{\infty }\left( \Omega
\times \left( 0,T\right) \right) $ and $y^{0}\in L^{2}\left( \Omega \right) $%
; we let $\omega $ be a nonempty open subset of $\Omega $; and we denote by $%
1_{\left\vert \cdot \right. }$ the characteristic function of a set
in the place where $\cdot $ stays.

Let $\tau \in \left[ 0,T\right) $. Let $E\subset \left( \tau
,T\right) $ be a measurable set of positive measure. Consider the
following parabolic
equation:%
\begin{equation}
\left\{
\begin{array}{ll}
\partial _{t}\psi -\Delta \psi +a\psi +B\cdot \nabla \psi =1_{\left\vert
\omega \times \left( \tau ,T\right) \right. }1_{\left\vert E\right.
}v &
\quad \text{in~}\Omega \times \left( 0,T\right) \ \text{,} \\
\psi =0 & \quad \text{on~}\partial \Omega \times \left( 0,T\right)
\text{ ,}
\\
\psi \left( \cdot ,0\right) =\psi ^{0} & \quad \text{in~}\Omega \ \text{,}%
\end{array}%
\right.  \tag{3.1}  \label{3.1}
\end{equation}%
where $v\in L^{\infty }\left( \Omega \times \left( 0,T\right) \right) $ and $%
\psi ^{0}\in L^{2}\left( \Omega \right) $. Then (\ref{3.1}) admits a
unique solution $\psi $ in $C\left( \left[ 0,T\right] ;L^{2}\left(
\Omega \right) \right) \cap L^{2}\left( 0,T;H_{0}^{1}\left( \Omega
\right) \right) $. The adjoint equation of (\ref{3.1}) is as:
\begin{equation}
\left\{
\begin{array}{ll}
-\partial _{t}\vartheta -\Delta \vartheta +\left(
a-\text{div}B\right) \vartheta -B\cdot \nabla \vartheta =0 & \quad
\text{in~}\Omega \times \left(
0,T\right) \ \text{,} \\
\vartheta =0 & \quad \text{on~}\partial \Omega \times \left(
0,T\right)
\text{ ,} \\
\vartheta \left( \cdot ,T\right) \in L^{2}\left( \Omega \right)
\text{ .} &
\end{array}%
\right.   \tag{3.2}  \label{3.2}
\end{equation}%
By Theorem 1.1, any solution $\vartheta $ to (\ref{3.2}) satisfies
\begin{equation}
\left\Vert \vartheta \left( \cdot ,0\right) \right\Vert
_{L^{2}\left( \Omega \right) }\leq \kappa \int_{\omega \times
E}\left\vert \vartheta \left( x,t\right) \right\vert dxdt\text{ ,}
\tag{3.3}  \label{3.3}
\end{equation}%
where the constant $\kappa $ is independent of $\vartheta $. This is
equivalent to the null-controllability from $E$: for any $\psi
^{0}\in L^{2}\left( \Omega \right) $, there is a $v\in L^{\infty
}\left( \Omega \times \left( 0,T\right) \right) $, with
\begin{equation}
\left\Vert v\right\Vert _{L^{\infty }\left( \Omega \times \left(
0,T\right) \right) }\leq \kappa \left\Vert \psi ^{0}\right\Vert
_{L^{2}\left( \Omega \right) }\text{ ,}  \tag{3.4}  \label{3.4}
\end{equation}%
such that the corresponding solution $\psi $ to (\ref{3.1})
satisfies $\psi \left( \cdot ,T\right) =0$ in $\Omega $. (See e.g.
\cite{W}). In general, such a $v$ is not unique.

\subsection{Norm optimal bang-bang control}

Consider the following parabolic equation:%
\begin{equation}
\left\{
\begin{array}{ll}
\partial _{t}y-\Delta y+ay+B\cdot \nabla y=1_{\left\vert \omega \times
\left( \tau ,T\right) \right. }f & \quad \text{in~}\Omega \times
\left(
0,T\right) \ \text{,} \\
y=0 & \quad \text{on~}\partial \Omega \times \left( 0,T\right) \text{ ,} \\
y\left( \cdot ,0\right) =y^{0} & \quad \text{in~}\Omega \ \text{,}%
\end{array}%
\right.  \tag{3.1.1}  \label{3.1.1}
\end{equation}%
where $f\in L^{\infty }\left( 0,T;L^{2}\left( \Omega \right) \right)
$. Then equation (\ref{3.1.1}) admits a unique solution $y$ in the
class of $C\left( \left[ 0,T \right] ;L^{2}\left( \Omega \right)
\right) \cap L^{2}\left( 0,T;H_{0}^{1}\left( \Omega \right) \right)
$. Write
\begin{equation*}
\mathcal{F}=\left\{ f\in L^{\infty }\left( 0,T;L^{2}\left( \Omega
\right) \right); y\left( \cdot ,T\right) =0\text{ in }\Omega
\right\} \text{ ,}
\end{equation*}
where $y$ is the solution of (\ref{3.1.1}) corresponding to $f$.

\begin{thm} \label{theorem3.1} There is a unique $f^{\ast }\in \mathcal{F}$\ such that
\begin{equation}
\left\Vert f^{\ast }\right\Vert _{L^{\infty }\left( \tau
,T;L^{2}\left( \Omega \right) \right) }=\underset{f\in
\mathcal{F}}{\min }\left\Vert
f\right\Vert _{L^{\infty }\left( \tau ,T;L^{2}\left( \Omega \right) \right) }%
\text{ .}  \tag{3.1.2}  \label{3.1.2}
\end{equation}%
Furthermore, $f^{\ast }$\ holds the bang-bang property:%
\begin{equation}
\left\Vert f^{\ast }\left( \cdot ,t\right) \right\Vert _{L^{2}\left(
\Omega \right) }=\left\Vert f^{\ast }\right\Vert _{L^{\infty }\left(
\tau ,T;L^{2}\left( \Omega \right) \right) }\text{ \quad\ for a.e.
}t\in \left( \tau ,T\right) \text{ .}  \tag{3.1.3}  \label{3.1.3}
\end{equation}
\end{thm}

\begin{rem} In the control theory of PDE, the equation (\ref{3.1.1}%
) is called a controlled system while $f$ is called a control. $f\in
\mathcal{F}$ means that the control $f$ in $L^{\infty }\left(
0,T;L^{2}\left( \Omega \right) \right) $ drives the solution $y$ of (\ref%
{3.1.1}) from $y^{0}$ to zero at time $T$. The property that
$\mathcal{F}$ is nonempty is called the null-controllability for
(\ref{3.1.1}). The quantity
\begin{equation}
\widetilde{M}=\underset{\widetilde{f}\in
\mathcal{F}}{\text{min}}\left\Vert \widetilde{f}\right\Vert
_{L^{\infty }\left( \tau ,T;L^{2}\left( \Omega \right) \right)
}\text{ }  \tag{3.1.4}  \label{3.1.4}
\end{equation}%
measures the best cost of such controls. The norm optimal control
problem (with respect to (\ref{3.1.1})) is to ask for a control
$f\in \mathcal{F}$ such that $\left\Vert f\right\Vert _{L^{\infty
}\left( \tau ,T;L^{2}\left( \Omega \right) \right) }=\widetilde{M}$.
Such a control is called a norm optimal control. The norm optimal
control problem has the bang-bang property if any norm optimal
control $f$ holds that $\left\Vert f\left( \cdot ,t\right)
\right\Vert _{L^{2}\left( \Omega \right) }=\widetilde{M}$ for a.e.
$t\in \left( \tau ,T\right) $. Theorem \ref{theorem3.1} presents
that the norm optimal problem has a unique optimal control and holds
the bang-bang property.
\end{rem}

\textbf{Proof of Theorem \ref{theorem3.1}.} We carry out the proof
by three steps as follows.

\textit{Step 1 .- Existence}. By the well-known result on the null
controllability of parabolic equations (see \cite{DFGZ}), we have that $%
\mathcal{F}\neq \emptyset $. Then by making use of the standard
argument of calculus of variations, we get the existence of such a
control $f\in \mathcal{F}$ satisfying $\left\Vert f\right\Vert
_{L^{\infty }\left( \tau ,T;L^{2}\left( \Omega \right) \right)
}=\widetilde{M}$.

\textit{Step 2 .- Bang-bang property}. We prove that if $f\in
\mathcal{F}$ satisfies (\ref{3.1.2}), then $f$ must hold
(\ref{3.1.3}). By seeking a contradiction, we suppose that
(\ref{3.1.3}) did not hold for some $f\in \mathcal{F}$ satisfying
(\ref{3.1.2}). Then there would be an $\varepsilon \in \left(
0,1\right) $ and a measurable set $E\subset \left( \tau ,T\right) $,
with a positive measure, such that
\begin{equation}
\left\Vert f\left( \cdot ,t\right) \right\Vert _{L^{2}\left( \Omega
\right) }\leq \widetilde{M}-\varepsilon \text{ \quad }\forall t\in
E\text{ .} \tag{3.1.5}  \label{3.1.5}
\end{equation}%
Here $\widetilde{M}$ is given by (\ref{3.1.4}). We claim that there are a $%
f_{\delta }\in L^{\infty }\left( 0,T;L^{2}\left( \Omega \right)
\right) $ with
\begin{equation}
\left\Vert f_{\delta }\right\Vert _{L^{\infty }\left( \tau
,T;L^{2}\left( \Omega \right) \right) }\leq \left( 1-\delta \right)
\widetilde{M}\text{ \quad for some }\delta \in \left( 0,1\right)
\text{ ,}  \tag{3.1.6} \label{3.1.6}
\end{equation}%
and a function $y_{\delta }$ with the property that
\begin{equation}
\left\{
\begin{array}{ll}
\partial _{t}y_{\delta }-\Delta y_{\delta }+ay_{\delta }+B\cdot \nabla
y_{\delta }=1_{\left\vert \omega \times \left( \tau ,T\right)
\right.
}f_{\delta } & \quad \text{in~}\Omega \times \left( 0,T\right) \ \text{,} \\
y_{\delta }=0 & \quad \text{on~}\partial \Omega \times \left(
0,T\right)
\text{ ,} \\
y_{\delta }\left( \cdot ,0\right) =y^{0} & \quad \text{in~}\Omega \
\text{,}
\\
y_{\delta }\left( \cdot ,T\right) =0 & \quad \text{in~}\Omega \ \text{.}%
\end{array}%
\right.  \tag{3.1.7}  \label{3.1.7}
\end{equation}%
The existence of such a triplet $\left( \delta, f_{\delta
},y_{\delta }\right) $ that satisfies (\ref{3.1.6}) and
(\ref{3.1.7}) clearly contradicts with the definition of
$\widetilde{M}$. Now, we prove the claim. Let $\delta \in \left(
0,1\right) $ (which will be determined later). By Theorem
\ref{theorem1.1} and its equivalence to the null-controllability
from $E$, there is a control $v_{\delta }\in L^{\infty }\left(
\Omega \times 0,T\right) $ such that the solution $\psi _{\delta }$
to
\begin{equation}
\left\{
\begin{array}{ll}
\partial _{t}\psi _{\delta }-\Delta \psi _{\delta }+a\psi _{\delta }+B\cdot
\nabla \psi _{\delta }=1_{\left\vert \omega \times \left( \tau
,T\right) \right. }1_{\left\vert E\right. }v_{\delta } & \quad
\text{in~}\Omega \times
\left( 0,T\right) \ \text{,} \\
\psi _{\delta }=0 & \quad \text{on~}\partial \Omega \times \left(
0,T\right)
\text{ ,} \\
\psi _{\delta }\left( \cdot ,0\right) =\delta y^{0} & \quad
\text{in~}\Omega
\ \text{,}%
\end{array}%
\right.  \tag{3.1.8}  \label{3.1.8}
\end{equation}%
satisfies $\psi _{\delta }\left( \cdot ,T\right) =0$ in $\Omega $.
Furthermore, there is a $\kappa >0$ (independent on $\delta $) such
that
\begin{equation}
\left\Vert v_{\delta }\right\Vert _{L^{\infty }\left(
0,T;L^{2}\left( \Omega \right) \right) }\leq \left\vert \Omega
\right\vert ^{1/2}\left\Vert v_{\delta }\right\Vert _{L^{\infty
}\left( \Omega \times \left( 0,T\right) \right) }\leq \kappa \delta
\left\Vert y^{0}\right\Vert _{L^{2}\left( \Omega \right) }\text{ .}
\tag{3.1.9}  \label{3.1.9}
\end{equation}%
Then we define $f_{\delta }$ by setting
\begin{equation}
f_{\delta }=\left( 1-\delta \right) f+1_{\left\vert E\right. }v_{\delta }%
\text{ .}  \tag{3.1.10}  \label{3.1.10}
\end{equation}%
By taking $\delta =\frac{\varepsilon }{\kappa \left\Vert
y^{0}\right\Vert _{L^{2}\left( \Omega \right) }+\varepsilon }$, one
can easily check that
\begin{equation}
\left\Vert f_{\delta }\left( \cdot ,t\right) \right\Vert
_{L^{2}\left( \Omega \right) }\leq \left( 1-\delta \right)
\widetilde{M}\text{ \quad\ for a.e. }t\in \left( \tau ,T\right)
\text{ .}  \tag{3.1.11}  \label{3.1.11}
\end{equation}%
On the other hand, one can verify that the function $\left( 1-\delta
\right)
y+\psi _{\delta }$ satisfies (\ref{3.1.7}). This, together with (\ref{3.1.11}%
), shows the claim.

\textit{Step 3 .- Uniqueness}. By the bang-bang property and the
parallelogram identity, we can easily check that the control $f\in
\mathcal{F}$ satisfying $\left\Vert f\right\Vert _{L^{\infty }\left(
\tau ,T;L^{2}\left( \Omega \right) \right) }=\widetilde{M}$ is
unique (see \cite[page 45]{F}).

This completes the proof.

\subsection{Time optimal bang-bang control}

Consider the following parabolic equation:%
\begin{equation}
\left\{
\begin{array}{ll}
\partial _{t}y-\Delta y+ay+B\cdot \nabla y=1_{\left\vert \omega \times
\left( \tau ,T\right) \right. }g & \quad \text{in~}\Omega \times
\left(
0,T\right) \ \text{,} \\
y=0 & \quad \text{on~}\partial \Omega \times \left( 0,T\right) \text{ ,} \\
y\left( \cdot ,0\right) =y^{0} & \quad \text{in~}\Omega \ \text{,}%
\end{array}%
\right.  \tag{3.2.1}  \label{3.2.1}
\end{equation}%
where $g\in L^{\infty }\left( 0,T;L^{2}\left( \Omega \right) \right)
$. Write
\begin{equation}
\mathcal{G}^{M}=\left\{ g\in L^{\infty }\left( 0,T;L^{2}\left(
\Omega \right) \right) ;\text{ }\left\Vert g\right\Vert _{L^{\infty
}\left( 0,T;L^{2}\left( \Omega \right) \right) }\leq M\right\}
\text{ ,}  \tag{3.2.2} \label{3.2.2}
\end{equation}%
where $M>0$. We define
\begin{equation}
\mathcal{P}^{M}=\left\{ \left( \tau ,g\right) \in \left[ 0,T\right)
\times \mathcal{G}^{M}; y\left( \cdot ,T\right) =0\text{ in }\Omega
\right\} \text{ ,} \tag{3.2.3}  \label{3.2.3}
\end{equation}
where $y$ is the solution of (\ref{3.2.1}) corresponding to $g$.

\begin{thm} \label{theorem3.3} Suppose that $\mathcal{P}^{M}\neq \emptyset $. If $\left( \tau
^{\ast },g^{\ast }\right) \in \mathcal{P}^{M}$\ is such that
\begin{equation}
\tau ^{\ast }\geq \tau \quad \ \text{for any pair }\left( \tau
,g\right) \in \mathcal{P}^{M}\text{ ,}  \tag{3.2.4}  \label{3.2.4}
\end{equation}%
then $g^{\ast }$\ holds the bang-bang property:%
\begin{equation}
\left\Vert g^{\ast }\left( \cdot ,t\right) \right\Vert _{L^{2}\left(
\Omega \right) }=M\text{ \quad\ for a.e. }t\in \left( \tau ^{\ast
},T\right) \text{ .}  \tag{3.2.5}  \label{3.2.5}
\end{equation}%
Furthermore, there is at most one such pairs $\left( \tau ^{\ast
},g^{\ast }\right) $.
\end{thm}

\begin{rem} It may happen that $\mathcal{P}^{M}=\emptyset
$. To guarantee that $\mathcal{P}^{M}\neq \emptyset $ for some
$T>0$, it is necessary to impose certain conditions on potentials
$a$ and $B$. For instance, it can be checked that one of the
following two conditions implies that $\mathcal{P}^{M}\neq \emptyset
$:

$\cdot$ $0\leq a-\frac{1}{2}$div$B+\lambda_{1}$ for a.e. $\left(
x,t\right) \in\Omega\times\left( 0,T\right) $ ;

$\cdot$ $\left\Vert a-\frac{1}{2}divB\right\Vert _{L^{\infty}\left(
\Omega\times\left( 0,T\right) \right) }\leq\lambda_{1}$ .

\noindent Here $\lambda_{1}>0$ denotes the first Dirichlet
eigenvalue.
\end{rem}

\begin{rem} There is a kind of time optimal control problem
whose aim is to delay initiation of active control (in a control
constraint set) as late as possible, such that the corresponding
solution (of a controlled system) reaches a target by a fixed ending
time (see e.g. \cite{MS}). In the current study, the controlled
system is (\ref{3.2.1}), where $g$ is viewed as a control; the
target is $\left\{ 0\right\} \subset L^{2}\left( \Omega \right) $;
the ending time is $T$; and the control constraint set is given by
(\ref{3.2.2}), where $M$ is regarded as a bound of controls. $\left(
\tau
,g\right) \in \mathcal{P}^{M}$ means that the control $g$ is not active in $%
\Omega \times \left( 0,\tau \right) $ and drives the solution of (\ref{3.2.1}%
) from $y^{0}$ to zero at time $T$. The time
\begin{equation*}
\tau ^{\ast }=\underset{\left( \tau ,g\right) \in \mathcal{P}^{M}}{\text{max}%
}\tau
\end{equation*}%
is called the optimal time; while a control $g^{\ast }$, with
$\left( \tau ^{\ast },g^{\ast }\right) \in \mathcal{P}^{M}$, is
called a time optimal control. Now from perspective of control
theory of PDE, Theorem \ref{theorem3.3} presents
that any time optimal control $g^{\ast }$ holds the bang-bang property: $%
\left\Vert g^{\ast }\left( \cdot ,t\right) \right\Vert _{L^{2}\left(
\Omega \right) }=M$ for a.e. $t\in \left( \tau ^{\ast },T\right) $.
It also shows that the optimal control, if it exists, is unique.
\end{rem}

\textbf{Proof of Theorem \ref{theorem3.3}.} The uniqueness of the
pair $\left( \tau
^{\ast },g^{\ast }\right) $ follows directly from the bang-bang property (%
\ref{3.2.5}) and the parallelogram identity (see \cite[page 45]{F}).
Thus, it remains to prove (\ref{3.2.5}). By contradiction, we
suppose that there was a pair $\left( \tau ^{\ast },g^{\ast }\right)
\in \mathcal{P}^{M}$ satisfying (\ref{3.2.4}) such that
(\ref{3.2.5}) did not hold. Then there
would be an $\varepsilon \in \left( 0,1\right) $ and a measurable set $%
\widetilde{E}\subset \left( \tau ^{\ast },T\right) $, with a
positive measure, such that
\begin{equation}
\left\Vert g^{\ast }\left( \cdot ,t\right) \right\Vert _{L^{2}\left(
\Omega \right) }\leq M-\varepsilon \text{ \quad }\forall t\in
\widetilde{E}\text{ .} \tag{3.2.6}  \label{3.2.6}
\end{equation}%
We claim that there are a $\delta \in \left( 0,1\right) $ and a pair
$\left( y,g\right) $ with $g\in \mathcal{G}^{M}$ such that
\begin{equation}
\left\{
\begin{array}{ll}
\partial _{t}y-\Delta y+ay+B\cdot \nabla y=1_{\left\vert \omega \times
\left( \tau ^{\ast }+\delta ,T\right) \right. }g & \quad
\text{in~}\Omega
\times \left( 0,T\right) \ \text{,} \\
y=0 & \quad \text{on~}\partial \Omega \times \left( 0,T\right) \text{ ,} \\
y\left( \cdot ,0\right) =y^{0} & \quad \text{in~}\Omega \ \text{,} \\
y\left( \cdot ,T\right) =0 & \quad \text{in~}\Omega \ \text{.}%
\end{array}%
\right.  \tag{3.2.7}  \label{3.2.7}
\end{equation}%
The existence of such a triplet $\left( \delta ,y,g\right) $ clearly
contradicts with (\ref{3.2.4}). To prove the claim, we first observe
that there is a $\delta _{0}\in \left( 0,1\right) $ such that the
measurable set
\begin{equation*}
E=\widetilde{E}\cap \left( \tau ^{\ast }+\delta _{0},T\right)
\end{equation*}%
has a positive measure. Then, it follows from (\ref{3.2.6}) that
\begin{equation}
\left\Vert g^{\ast }\left( \cdot ,t\right) \right\Vert _{L^{2}\left(
\Omega \right) }\leq M-\varepsilon \text{ \quad }\forall t\in
E\text{ .} \tag{3.2.8}  \label{3.2.8}
\end{equation}%
Let $\delta \in \left( 0,\delta _{0}\right) $, which will be
determined later. By solving the equation:
\begin{equation}
\left\{
\begin{array}{ll}
\partial _{t}z-\Delta z+az+B\cdot \nabla z=-1_{\left\vert \omega \times
\left( \tau ^{\ast },\tau ^{\ast }+\delta \right) \right. }g^{\ast }
& \quad
\text{in~}\Omega \times \left( 0,\tau ^{\ast }+\delta \right) \ \text{,} \\
z=0 & \quad \text{on~}\partial \Omega \times \left( 0,\tau ^{\ast
}+\delta
\right) \text{ ,} \\
z\left( \cdot ,0\right) =0 & \quad \text{in~}\Omega \ \text{,}%
\end{array}%
\right.  \tag{3.2.9}  \label{3.2.9}
\end{equation}%
we get that
\begin{equation}
\begin{array}{ll}
\left\Vert z\left( \cdot ,\tau ^{\ast }+\delta \right) \right\Vert
_{L^{2}\left( \Omega \right) } & \leq c_{0}\left\Vert g^{\ast
}\right\Vert _{L^{1}\left( \tau ^{\ast },\tau ^{\ast }+\delta
;L^{2}\left( \Omega \right)
\right) } \\
& \leq c_{0}M\delta \text{ ,}%
\end{array}
\tag{3.2.10}  \label{3.2.10}
\end{equation}%
where $c_{0}>0$ is independent on $\delta $.

Next, by Theorem \ref{theorem1.1} and its equivalence to the
null-controllability from $E$, there is a control $v\in L^{\infty
}\left( \Omega \times \left( \tau ^{\ast }+\delta ,T\right) \right)
$ such that the solution $\psi $ to the equation:
\begin{equation}
\left\{
\begin{array}{ll}
\partial _{t}\psi -\Delta \psi +a\psi +B\cdot \nabla \psi =1_{\left\vert
\omega \times \left( \tau ^{\ast }+\delta ,T\right) \right.
}1_{\left\vert E\right. }v & \quad \text{in~}\Omega \times \left(
\tau ^{\ast
}+\delta ,T\right) \ \text{,} \\
\psi =0 & \quad \text{on~}\partial \Omega \times \left( \tau ^{\ast
}+\delta
,T\right) \text{ ,} \\
\psi \left( \cdot ,\tau ^{\ast }+\delta \right) =z\left( \cdot ,\tau
^{\ast
}+\delta \right)  & \quad \text{in~}\Omega \ \text{,}%
\end{array}%
\right.   \tag{3.2.11}  \label{3.2.11}
\end{equation}%
satisfies $\psi \left( \cdot ,T\right) =0$ in $\Omega $.
Furthermore, it holds that
\begin{equation}
\left\Vert v\right\Vert _{L^{\infty }\left( \tau ^{\ast }+\delta
,T;L^{2}\left( \Omega \right) \right) }\leq \left\vert \Omega
\right\vert ^{1/2}\left\Vert v\right\Vert _{L^{\infty }\left( \Omega
\times \left( \tau ^{\ast }+\delta ,T\right) \right) }\leq \kappa
\left\Vert z\left( \cdot
,\tau ^{\ast }+\delta \right) \right\Vert _{L^{2}\left( \Omega \right) }%
\text{ ,}  \tag{3.2.12}  \label{3.2.12}
\end{equation}%
for some $\kappa >0$ independent on $\delta $. Combining the above
estimate with (\ref{3.2.10}), we can find a constant $c>0$,
independent on $\delta $, such that
\begin{equation}
\left\Vert v\right\Vert _{L^{\infty }\left( \tau ^{\ast }+\delta
,T;L^{2}\left( \Omega \right) \right) }\leq c\delta \text{ .}
\tag{3.2.13} \label{3.2.13}
\end{equation}%
Now, we define%
\begin{equation}
w\left( \cdot ,t\right) =\left\vert
\begin{array}{ll}
z\left( \cdot ,t\right)  & \text{if }t\in \left[ 0,\tau ^{\ast }+\delta %
\right] \text{ ,} \\
\psi \left( \cdot ,t\right)  & \text{if }t\in \left( \tau ^{\ast }+\delta ,T%
\right] \text{ .}%
\end{array}%
\right.   \tag{3.2.14}  \label{3.2.14}
\end{equation}%
Clearly, $w\left( \cdot ,0\right) =0$ and $w\left( \cdot ,T\right) =0$ in $%
\Omega $. Let $y^{\ast }$ be the solution of (\ref{3.2.1}) with
$\left( \tau
,g\right) =\left( \tau ^{\ast },g^{\ast }\right) $. Thus it holds that $%
y^{\ast }\left( \cdot ,T\right) =0$ in $\Omega $. Further, one can
easily check that the function $y^{\ast }+w$ solves (\ref{3.2.7})
with
\begin{equation}
g\left( \cdot ,t\right) =\left\vert
\begin{array}{ll}
0 & \text{if }t\in \left[ 0,\tau ^{\ast }+\delta \right] \text{ ,} \\
g^{\ast }\left( \cdot ,t\right) +1_{\left\vert E\right. }v\left(
\cdot
,t\right)  & \text{if }t\in \left( \tau ^{\ast }+\delta ,T\right] \text{ .}%
\end{array}%
\right.   \tag{3.2.15}  \label{3.2.15}
\end{equation}%
Finally, we take $\delta \in \left( 0,\delta _{0}\right) $ such that $%
c\delta \leq \varepsilon $. Then it holds that $g\in
\mathcal{G}^{M}$. Indeed, it follows from (\ref{3.2.15}),
(\ref{3.2.8}) and (\ref{3.2.13}) that
\begin{equation*}
\begin{array}{ll}
\left\Vert g\left( \cdot ,t\right) \right\Vert _{L^{2}\left( \Omega
\right) } & \leq \left\Vert g^{\ast }\left( \cdot ,t\right)
\right\Vert _{L^{2}\left( \Omega \right) }+\left\Vert 1_{\left\vert
E\right. }v\left( \cdot ,t\right)
\right\Vert _{L^{2}\left( \Omega \right) } \\
& \leq \left\vert
\begin{array}{ll}
M-\varepsilon +c\delta & \text{a.e. if }t\in E\cap \left( \tau
^{\ast
}+\delta ,T\right) \\
M & \text{a.e. if }t\notin E\cap \left( \tau ^{\ast }+\delta ,T\right) \\
0 & \text{if }t\in \left( 0,\tau ^{\ast }+\delta \right)%
\end{array}%
\right. \\
& \leq M\quad \text{for a.e. }t\in \left( 0,T\right) \text{ .}%
\end{array}%
\end{equation*}%
This completes the proof.

\section*{Appendix}

\subsection*{Proof of Proposition \ref{proposition2.1}}

Since $\left\vert E\right\vert >0$, almost every point of $E$ is a
point of density of $E\subset \left( 0,T\right) $. Let $\ell \in
\left( 0,T\right) $ be such a point. Then it holds that
\begin{equation}
\frac{\left\vert E^{c}\cap \left( \ell -\theta ,\ell +\theta \right)
\right\vert }{\left\vert \left( \ell -\theta ,\ell +\theta \right)
\right\vert }\rightarrow 0\text{ and }\frac{\left\vert E\cap \left(
\ell -\theta ,\ell +\theta \right) \right\vert }{\left\vert \left(
\ell -\theta ,\ell +\theta \right) \right\vert }\rightarrow 1\text{
as }\theta \rightarrow 0\text{ .}  \tag{A.1.1}  \label{A.1.1}
\end{equation}%
Let $z>1$. Let $0<\epsilon \leq $min$\left( \frac{z-1}{1+3z},\frac{1}{3}%
\right) $ which implies that
\begin{equation}
\frac{\epsilon }{1-\epsilon }\left( \frac{1+z}{z-1}\right) \leq \frac{1}{2}%
\text{ and }\left( 1+\frac{\epsilon }{1-\epsilon }\right) \leq \frac{3}{2}%
\text{ .}  \tag{A.1.2}  \label{A.1.2}
\end{equation}%
Then by (\ref{A.1.1}), there exists $\theta _{o}=\theta _{o}\left(
\epsilon \right) >0$ such that for any $\theta <\theta _{o}$,
\begin{equation*}
\frac{\left\vert E^{c}\cap \left( \ell -\theta ,\ell +\theta \right)
\right\vert }{\left\vert \left( \ell -\theta ,\ell +\theta \right)
\right\vert }<\epsilon \text{ and }1-\epsilon <\frac{\left\vert
E\cap \left( \ell -\theta ,\ell +\theta \right) \right\vert
}{\left\vert \left( \ell -\theta ,\ell +\theta \right) \right\vert
}\text{ ,}
\end{equation*}%
which imply that
\begin{equation}
\left\vert E^{c}\cap \left( \ell -\theta ,\ell +\theta \right) \right\vert <%
\frac{\epsilon }{1-\epsilon }\left\vert E\cap \left( \ell -\theta
,\ell +\theta \right) \right\vert \text{ .}  \tag{A.1.3}
\label{A.1.3}
\end{equation}%
Write $\widetilde{\theta _{o}}=$min$\left( \theta _{o},T-\ell
\right) $. Let $\ell _{1}$ be such that $\ell <\ell _{1}<\ell
+\widetilde{\theta _{o}}\leq T $. Define $\left\{ \ell _{m}\right\}
_{m\geq 1}$ by (\ref{2.1.1}). Clearly,
\begin{equation}
\ell _{m}-\ell <\ell _{m-1}-\ell <\cdots <\ell _{2}-\ell <\ell _{1}-\ell <%
\widetilde{\theta _{o}}\leq \theta _{o}  \tag{A.1.4}  \label{A.1.4}
\end{equation}%
and
\begin{equation}
\ell _{m+1}-\ell _{m+2}=\frac{1}{z^{m+1}}\left( z-1\right) \left(
\ell _{1}-\ell \right) \text{ .}  \tag{A.1.5}  \label{A.1.5}
\end{equation}%
Then%
\begin{equation}
\begin{array}{ll}
\ell _{m}-\ell _{m+1} & =\left\vert E^{c}\cap \left( \ell
_{m+1},\ell _{m}\right) \right\vert +\left\vert E\cap \left( \ell
_{m+1},\ell
_{m}\right) \right\vert \\
& \leq \left\vert E^{c}\cap \left( 2\ell -\ell _{m},\ell _{m}\right)
\right\vert +\left\vert E\cap \left( \ell _{m+1},\ell _{m}\right)
\right\vert
\\
& \leq \frac{\epsilon }{1-\epsilon }\left\vert E\cap \left( 2\ell
-\ell _{m},\ell _{m}\right) \right\vert +\left\vert E\cap \left(
\ell _{m+1},\ell
_{m}\right) \right\vert \text{ .}%
\end{array}
\tag{A.1.6}  \label{A.1.6}
\end{equation}%
The first inequality in (\ref{A.1.6}) follows from (\ref{2.1.1});
while in the second inequality of (\ref{A.1.6}), we used
(\ref{A.1.3}), with $\theta =\ell _{m}-\ell $, and (\ref{A.1.4}).
Thus we have that
\begin{equation}
\begin{array}{ll}
\ell _{m}-\ell _{m+1} & \leq \frac{\epsilon }{1-\epsilon }\left[
\left\vert E\cap \left( 2\ell -\ell _{m},\ell _{m+1}\right)
\right\vert +\left\vert E\cap \left( \ell _{m+1},\ell _{m}\right)
\right\vert \right] +\left\vert
E\cap \left( \ell _{m+1},\ell _{m}\right) \right\vert \\
& \leq \left( 1+\frac{\epsilon }{1-\epsilon }\right) \left\vert
E\cap \left( \ell _{m+1},\ell _{m}\right) \right\vert
+\frac{\epsilon }{1-\epsilon }\left[
\ell _{m+1}-\left( 2\ell -\ell _{m}\right) \right] \text{ .}%
\end{array}
\tag{A.1.7}  \label{A.1.7}
\end{equation}%
Besides, it follows from (\ref{2.1.1}) and (\ref{A.1.5}) that
\begin{equation*}
\begin{array}{ll}
\ell _{m+1}-\left( 2\ell -\ell _{m}\right) & =\frac{1}{z^{m}}\left(
1+z\right) \left( \ell _{1}-\ell \right) \\
& =\frac{1+z}{z-1}\left( \ell _{m}-\ell _{m+1}\right) \text{ .}%
\end{array}%
\end{equation*}%
This, along with (\ref{A.1.7}) and (\ref{A.1.2}), leads to
(\ref{2.1.2}).

\subsection*{Some inequalities}

Suppose $a\in L^{\infty }\left( 0,T;L^{q}\left( \Omega \right)
\right) $ where $q\geq 2$ for $n=1$, and $q>n$ for $n\geq 2$. Let
\begin{equation*}
p=\left\vert
\begin{array}{ll}
\frac{2n}{q} & \text{if }n<q\leq 2n \\
1 & \text{if }2n\leq q%
\end{array}%
\right.
\end{equation*}%
Then, for each $\varepsilon >0$, there is $C_{\left( \varepsilon
,\Omega ,n,q\right) }>0$ such that for any $\phi \in H_{0}^{1}\left(
\Omega \right) $ and for a.e. $t\in \left( 0,L\right) \subset \left[
0,T\right] $,
\begin{equation}
\begin{array}{ll}
& \quad \displaystyle\int_{\Omega }\left\vert a\left( x,t\right)
\phi \left(
x\right) \right\vert ^{2}G_{\lambda }\left( x,t\right) dx \\
& \leq \varepsilon \int_{\Omega }\left\vert \nabla \phi \left(
x\right)
\right\vert ^{2}G_{\lambda }\left( x,t\right) dx \\
& \quad +C_{\left( \varepsilon ,\Omega ,n,q\right) }\left(
\left\Vert a\right\Vert _{L^{\infty }\left( 0,L;L^{q}\left( \Omega
\right) \right) }^{4/\left( 2-p\right) }+\frac{\left\Vert
a\right\Vert _{L^{\infty }\left(
0,L;L^{q}\left( \Omega \right) \right) }^{2}}{\left( L-t+\lambda \right) ^{p}%
}\right) \int_{\Omega }\left\vert \phi \left( x\right) \right\vert
^{2}G_{\lambda }\left( x,t\right) dx%
\end{array}
\tag{A.2.1}  \label{A.2.1}
\end{equation}%
and%
\begin{equation}
\begin{array}{ll}
& \quad \displaystyle\int_{\Omega }a\left( x,t\right) \left\vert
\phi \left(
x\right) \right\vert ^{2}G_{\lambda }\left( x,t\right) dx \\
& \leq \varepsilon \int_{\Omega }\left\vert \nabla \phi \left(
x\right)
\right\vert ^{2}G_{\lambda }\left( x,t\right) dx \\
& \quad +C_{\left( \varepsilon ,\Omega ,n,q\right) }\left(
\left\Vert a\right\Vert _{L^{\infty }\left( 0,L;L^{q}\left( \Omega
\right) \right) }^{2}+\frac{\left\Vert a\right\Vert _{L^{\infty
}\left( 0,L;L^{q}\left( \Omega \right) \right) }}{L-t+\lambda
}\right) \int_{\Omega }\left\vert \phi
\left( x\right) \right\vert ^{2}G_{\lambda }\left( x,t\right) dx\text{ .}%
\end{array}
\tag{A.2.2}  \label{A.2.2}
\end{equation}

Proof of (\ref{A.2.1}) and (\ref{A.2.2}). Notice that $1\leq p<2$.
In the case where $n\geq 2$, it holds that
\begin{equation*}
\begin{array}{ll}
& \quad \displaystyle\int_{\Omega }\left\vert a\phi \right\vert
^{2}G_{\lambda }dx \\
& \leq \left\Vert a^{2}\right\Vert _{L^{\frac{n}{p}}\left( \Omega
\right) }\left\Vert \phi ^{2}G_{\lambda }\right\Vert
_{L^{\frac{n}{n-p}}\left(
\Omega \right) }\text{\quad (by H\"{o}lder inequality)} \\
& \leq \left\Vert a\right\Vert _{L^{\frac{2n}{p}}\left( \Omega
\right)
}^{2}\left\Vert \left( \phi ^{2}G_{\lambda }\right) ^{\frac{1}{p}%
}\right\Vert _{L^{\frac{pn}{n-p}}\left( \Omega \right) }^{p} \\
& \leq C_{\left( \Omega ,n,q\right) }\left\Vert a\right\Vert
_{L^{q}\left( \Omega \right) }^{2}\left\Vert \nabla \left( \left(
\phi ^{2}G_{\lambda
}\right) ^{\frac{1}{p}}\right) \right\Vert _{L^{p}\left( \Omega \right) }^{p}%
\text{\quad (by Sobolev inequality)} \\
& \leq C_{\left( \Omega ,n,q\right) }\left\Vert a\right\Vert
_{L^{q}\left( \Omega \right) }^{2}\displaystyle\int_{\Omega }\left(
\left\vert \phi \right\vert ^{2-p}\left\vert \nabla \phi \right\vert
^{p}G_{\lambda }+\left\vert \phi \right\vert ^{2}\left( G_{\lambda
}\right)
^{1-p}\left\vert \nabla G_{\lambda }\right\vert ^{p}\right) dx \\
& \leq C_{\left( \Omega ,n,q\right) }\left\Vert a\right\Vert
_{L^{q}\left( \Omega \right) }^{2}\left( \displaystyle\int_{\Omega
}\left\vert \phi
\right\vert ^{2}G_{\lambda }dx\right) ^{\frac{2-p}{2}}\left( \displaystyle%
\int_{\Omega }\left\vert \nabla \phi \right\vert ^{2}G_{\lambda }dx\right) ^{%
\frac{p}{2}}\text{ (by H\"{o}lder inequality)} \\
& \quad +C_{\left( \Omega ,n,q\right) }\left\Vert a\right\Vert
_{L^{q}\left(
\Omega \right) }^{2}\left( \frac{1}{\left\vert L-t+\lambda \right\vert ^{p}}%
\displaystyle\int_{\Omega }\left\vert \phi \right\vert
^{2}G_{\lambda
}dx\right)  \\
& \leq C_{\left( \Omega ,n,q\right) }\left\Vert a\right\Vert
_{L^{q}\left( \Omega \right) }^{2}\left( \varepsilon
\displaystyle\int_{\Omega }\left\vert
\nabla \phi \right\vert ^{2}G_{\lambda }dx+\left( \frac{1}{\varepsilon ^{%
\frac{p}{2-p}}}+\frac{1}{\left\vert L-t+\lambda \right\vert ^{p}}\right) %
\displaystyle\int_{\Omega }\left\vert \phi \right\vert
^{2}G_{\lambda
}dx\right) \text{ ,}%
\end{array}%
\end{equation*}%
and%
\begin{equation*}
\begin{array}{ll}
& \quad \displaystyle\int_{\Omega }a\left\vert \phi \right\vert
^{2}G_{\lambda }dx \\
& \leq \left\Vert a\right\Vert _{L^{n}\left( \Omega \right)
}\left\Vert \phi
^{2}G_{\lambda }\right\Vert _{L^{\frac{n}{n-1}}\left( \Omega \right) }\text{%
\quad (by H\"{o}lder inequality)} \\
& \leq C_{\left( \Omega ,n\right) }\left\Vert a\right\Vert
_{L^{n}\left( \Omega \right) }\left\Vert \nabla \left( \phi
^{2}G_{\lambda }\right)
\right\Vert _{L^{1}\left( \Omega \right) }\text{\quad (by Sobolev inequality)%
} \\
& \leq C_{\left( \Omega ,n,q\right) }\left\Vert a\right\Vert
_{L^{q}\left( \Omega \right) }\displaystyle\int_{\Omega }\left(
\left\vert \phi \right\vert \left\vert \nabla \phi \right\vert
G_{\lambda }+\left\vert \phi
\right\vert ^{2}\left\vert \nabla G_{\lambda }\right\vert \right) dx\text{%
\quad (by H\"{o}lder inequality)} \\
& \leq C_{\left( \Omega ,n,q\right) }\left\Vert a\right\Vert
_{L^{q}\left( \Omega \right) }\left( \varepsilon
\displaystyle\int_{\Omega }\left\vert
\nabla \phi \right\vert ^{2}G_{\lambda }dx+\left( \frac{1}{\varepsilon }+%
\frac{1}{\left\vert L-t+\lambda \right\vert }\right) \displaystyle%
\int_{\Omega }\left\vert \phi \right\vert ^{2}G_{\lambda }dx\right) \text{.}%
\end{array}%
\end{equation*}%
In the case when $n=1$, it stands that
\begin{equation*}
\begin{array}{ll}
& \quad \displaystyle\int_{\Omega }\left\vert a\phi \right\vert
^{2}G_{\lambda }dx \\
& \leq \left\Vert a^{2}\right\Vert _{L^{1}\left( \Omega \right)
}\left\Vert
\phi ^{2}G_{\lambda }\right\Vert _{L^{\infty }\left( \Omega \right) }\text{%
\quad (by H\"{o}lder inequality)} \\
& \leq C_{\left( \Omega \right) }\left\Vert a\right\Vert
_{L^{2}\left( \Omega \right) }^{2}\left\Vert \nabla \left( \phi
^{2}G_{\lambda }\right)
\right\Vert _{L^{1}\left( \Omega \right) }\text{\quad (by Sobolev inequality)%
} \\
& \leq C_{\left( \Omega \right) }\left\Vert a\right\Vert
_{L^{2}\left( \Omega \right) }^{2}\displaystyle\int_{\Omega }\left(
\left\vert \phi \right\vert \left\vert \nabla \phi \right\vert
G_{\lambda }+\left\vert \phi
\right\vert ^{2}\left\vert \nabla G_{\lambda }\right\vert \right) dx \\
& \leq C_{\left( \Omega \right) }\left\Vert a\right\Vert
_{L^{2}\left( \Omega \right) }^{2}\left( \int_{\Omega }\left\vert
\phi \right\vert ^{2}G_{\lambda }dx\right) ^{1/2}\left( \int_{\Omega
}\left\vert \nabla \phi \right\vert ^{2}G_{\lambda }dx\right)
^{1/2}\text{\quad (by H\"{o}lder
inequality)} \\
& \quad +C_{\left( \Omega \right) }\left\Vert a\right\Vert
_{L^{2}\left(
\Omega \right) }^{2}\left( \frac{1}{\left\vert L-t+\lambda \right\vert }%
\displaystyle\int_{\Omega }\left\vert \phi \right\vert
^{2}G_{\lambda
}dx\right)  \\
& \leq C_{\left( \Omega \right) }\left\Vert a\right\Vert
_{L^{2}\left( \Omega \right) }^{2}\left( \varepsilon
\displaystyle\int_{\Omega }\left\vert
\nabla \phi \right\vert ^{2}G_{\lambda }dx+\left( \frac{1}{\varepsilon }+%
\frac{1}{\left\vert L-t+\lambda \right\vert }\right) \displaystyle%
\int_{\Omega }\left\vert \phi \right\vert ^{2}G_{\lambda }dx\right) \text{ ,}%
\end{array}%
\end{equation*}%
and%
\begin{equation*}
\begin{array}{ll}
& \quad \displaystyle\int_{\Omega }a\left\vert \phi \right\vert
^{2}G_{\lambda }dx \\
& \leq \left\Vert a\right\Vert _{L^{1}\left( \Omega \right)
}\left\Vert \phi ^{2}G_{\lambda }\right\Vert _{L^{\infty }\left(
\Omega \right) }\text{\quad
(by H\"{o}lder inequality)} \\
& \leq C_{\left( \Omega \right) }\left\Vert a\right\Vert
_{L^{1}\left( \Omega \right) }\left\Vert \nabla \left( \phi
^{2}G_{\lambda }\right)
\right\Vert _{L^{1}\left( \Omega \right) }\text{\quad (by Sobolev inequality)%
} \\
& \leq C_{\left( \Omega \right) }\left\Vert a\right\Vert
_{L^{1}\left( \Omega \right) }\displaystyle\int_{\Omega }\left(
\left\vert \phi \right\vert \left\vert \nabla \phi \right\vert
G_{\lambda }+\left\vert \phi
\right\vert ^{2}\left\vert \nabla G_{\lambda }\right\vert \right) dx\text{%
\quad (by H\"{o}lder inequality)} \\
& \leq C_{\left( \Omega \right) }\left\Vert a\right\Vert
_{L^{2}\left( \Omega \right) }\left( \varepsilon
\displaystyle\int_{\Omega }\left\vert
\nabla \phi \right\vert ^{2}G_{\lambda }dx+\left( \frac{1}{\varepsilon }+%
\frac{1}{\left\vert L-t+\lambda \right\vert }\right) \displaystyle%
\int_{\Omega }\left\vert \phi \right\vert ^{2}G_{\lambda }dx\right)
\text{ }\forall \varepsilon >0\text{ .}
\end{array}%
\end{equation*}

\bigskip
\footnotesize \noindent\textit{Acknowledgments.} This work was
partially supported by the National Basis Research Program of China
(973 Program) under grant 2011CB808002 and the National Science
Foundation of China under grant 11161130003.


\begin{thebibliography}{SK}




\normalsize
\baselineskip=17pt


\bibitem[BT]{BT} Bardos, C., Tartar, L.: Sur l'unicit\'{e} retrograde
des \'{e}quations paraboliques et quelques questions voisines, Arch.
Rational Mech. Anal. \textbf{50}, 10--25 (1973)

\bibitem[DFGZ]{DFGZ} Doubova, A., Fernandez-Cara, E., Gonzalez-Burgos, M.,
Zuazua, E.: On the controllability of parabolic systems with a
nonlinear term involving the state and the gradient, SIAM J. Control
Optim. \textbf{41}, 798--819 (2002)

\bibitem[DZZ]{DZZ} Duyckaerts, T., Zhang, X., Zuazua, E.: On the
optimality of the observability inequalities for parabolic and
hyperbolic systems with potentials, Ann. Inst. H. Poincar\'{e} (C)
Non Linear Analysis \textbf{25}, 1--41 (2008)

\bibitem[EFV]{EFV} Escauriaza, L., Fernandez, F.J., Vessella, S.:
Doubling properties of caloric functions, Applicable analysis
\textbf{85} (1-2), 205--223 (2006)

\bibitem[F]{F} Fattorini, H.: Infinite Dimensional Linear Control Systems;
the Time Optimal and Norm Optimal Problems, North-Holland
Mathematical Studies vol 201, Elsevier, Amsterdam (2005)

\bibitem[FI]{FI} Fursikov, A.V., Imanuvilov, O.Y.: Controllability of
Evolution Equations. Lecture Notes Series, 34. Seoul National
University, Research Institute of Mathematics, Global Analysis
Research Center, Seoul (1996)

\bibitem[K]{K} Kenig, C.: Quantitative unique continuation,
logarithmic convexity of gaussian means and Hardy's uncertainty
principle, Proc. of Symp. of Pure Math., Volume in celebration of V.
Mazya's 70 birthday, 79, 207--227 (2008)

\bibitem[KT]{KT} Koch, H., Tataru, D.: Carleman estimates and unique
continuation for second order parabolic equations with non smooth
coefficients, Comm. Partial\ Differential Equations \textbf{34},
305--366 (2009)

\bibitem[LR]{LR} Lebeau, G., Robbiano, L.: Contr\^{o}le exact de l'
\'{e}quation de la chaleur, Comm. Partial\ Differential Equations
\textbf{20}, 335--356 (1995)

\bibitem[L]{L} Lin, F.H.: Remarks on a backward parabolic problem,
Methods and Applications of Analysis \textbf{10} (2), 245-252 (2003)

\bibitem[M]{M} Miller, L.: A direct Lebeau-Robbiano strategy for the
observability of heat-like semigroups, Discrete and Continuous
Dynamical Systems Series B \textbf{14} (4), 1465--1485 (2010)

\bibitem[MS]{MS} Mizel, V., Seidman, T.: An abstract 'bang-bang
principle' and time optimal boundary control of the heat equation,
SIAM J. Control Optim. \textbf{35}, 1204--1216 (1997)

\bibitem[PW]{PW} Phung, K.D., Wang, G.: Quantitative unique
continuation for the semilinear heat equation in a convex domain,
Journal of Functional Analysis \textbf{259} (5), 1230--1247 (2010)

\bibitem[P]{P} Poon, C.C.: Unique continuation for parabolic equations
, Comm. Partial\ Differential Equations \textbf{21}, 521--539 (1996)

\bibitem[W]{W} Wang, G.: L$^{\infty}$-null controllability
for the heat equation and its consequences for the time optimal
control problem, SIAM J. Control Optim. \textbf{47}, 1701--1720
(2008)

\end{thebibliography}
\end{document}